\documentclass{article}
\usepackage[numbers]{natbib}

\usepackage{framed,multirow}

\usepackage{amssymb}
\usepackage{latexsym}

\usepackage{url}

\usepackage{lmodern}
\usepackage[T1]{fontenc}
\usepackage{setspace}

\usepackage{pgfplots}
\pgfplotsset{compat=1.11}
\usepackage{tikz,tkz-euclide}
\usepackage{tikz-cd}
\usepackage{color}
\usepackage{tikz}
\usetikzlibrary{arrows,calc,patterns}
\usepackage{pgfplotstable}

\usepackage{rotating}
\usepackage{amsmath,mathtools,bm,cancel,relsize}
\usepackage{graphics,psfrag}
\usepackage{subcaption} 
\usepackage{empheq}
\usepackage{colortbl}
\usepackage{hyperref}
\usepackage{fancyref}
\usepackage{enumitem}
\usepackage{amsthm}
\usepackage[capitalize]{cleveref}

\usepackage{authblk}

\newtheorem{thm}{Theorem}[section]

\newtheorem{remark}[thm]{Remark}

\numberwithin{equation}{section}

\newcommand{\R}{\mathbf{R}}  
\newcommand{\N}{\mathbf{N}}  

\DeclareMathOperator\spa{span} 

\newcommand{\pf}{\iota}

\newcommand{\Xgradp}[1]{X_P^{\grad}({#1})}
\newcommand{\Xcurlp}[1]{X_P^{\curl}({#1})}
\newcommand{\Xdivp}[1]{X_P^{\div}({#1})}
\newcommand{\Xgradd}[1]{\tilde X_P^{\grad}({#1})}
\newcommand{\Xcurld}[1]{\tilde X_P^{\curl}({#1})}
\newcommand{\Xdivd}[1]{\tilde X_P^{\div}({#1})}

\def\triangulation{\mathcal{T}_h}
\def\dualmesh{{\tilde{\mathcal{T}}}_h}
\def\quadmesh{\mathcal{K}_h}

\def\sumoverallfacetslocald{\sum_{F \in \partial \tilde T}}
\def\sumoverallfacetslocalp{\sum_{F \in \partial T}}
\def\sumoverallptrigs{\sum_{K \in \triangulation}}
\def\sumoveralldtrigs{\sum_{\tilde{T} \in \dualmesh}}

\def\sumoverallelementsindtrig{\sum_{K \subset \tilde T}}
\def\sumoverallelementsinptrig{\sum_{K \subset T}}
\def\normal{\hat{\bm n}}

\def\grad{\operatorname{grad}}                                
\def\curl{\operatorname{curl}}                                
\def\div{\operatorname{div}}                                
\def\rot{\operatorname{rot}}                                

\title{Mass lumping the dual cell method to arbitrary polynomial degree for acoustic and electromagnetic waves}
\author[1]{Markus Wess}
\author[2]{Bernard Kapidani}
\author[3]{Lorenzo Codecasa}
\author[1]{Joachim Schoberl}

\affil[1]{Institute of Analysis and Scientific Computing, Technische Universit\unexpanded{\"a}t Wien, A-1040, Vienna, Austria.\authorcr
    \tt markus.wess@tuwien.ac.at, joachim.schoeberl@tuwien.ac.at }
\affil[2]{Institute of Mathematics,  \unexpanded{\'{E}}cole Polytechnique F\unexpanded{\'{e}}d\unexpanded{\'{e}}rale Lausanne, CH-1015 Lausanne, Switzerland \authorcr
        \tt bernard.kapidani@epfl.ch }
\affil[3]{Dipartimento di Elettronica, Informatica e Bioingegneria, Politecnico di Milano, I-20133 Milano, Italy \authorcr
    \tt lorenzo.codecasa@polimi.it }

\date{\today}

\begin{document}
\maketitle
\begin{abstract}
  We present a fundamental improvement of a high polynomial degree time domain cell method recently introduced by the last three authors. The published work introduced a method featuring block-diagonal system matrices where the block size and conditioning scaled poorly with respect to polynomial degree. The issue is herein bypassed by the construction of new basis functions exploiting quadrature rule based mass lumping techniques for arbitrary polynomial degrees in two dimensions for the Maxwell equations and the acoustic wave equation in the first order velocity pressure formulation. We characterize the degrees of freedom of all new discrete approximation spaces we employ for differential forms and show that the resulting block diagonal (inverse) mass matrices have block sizes independent of the polynomial degree. We demonstrate on an extensive number of examples how the new technique is applicable and efficient for large scale computations. 

\end{abstract}

\section{Introduction}\label{sec:intro}
When solving time dependent initial boundary value problems for hyperbolic partial differential equations such as the Maxwell system or the acoustic wave equation, the most used choice for the space discretisation is finite differences (usually the second order accurate version, on staggered grids).
This is due to their massively parallelisable nature and the fact that they are easy to implement. Especially the first point is of great importance for the simulation of large scale problems, e.g., in the context of seismic imaging or the simulation of electromagnetic waves in photonic crystals\citep{joannopoulosPhotonicCrystalsMolding2008,tafloveAdvancesFDTDComputational2013}.
Nevertheless, since the work of Hesthaven \& Warburton \citep{hesthavenNodalDiscontinuousGalerkin2008} on discontinuous Galerkin (DG) Finite Element Methods (FEM) there has been a revitalized enthusiasm in using variational methods to discretise the Maxwell system in a way that leads to block diagonal mass matrices even on unstructured grids.

The present manuscript fits in this framework and, while building originally on low order Finite Integration Techniques (FIT, \citep{weilandTimeDomainElectromagnetic1996,codecasaExplicitConsistentConditionally2008}), is a instead high-order accurate like the DG approach and acts as a follow up on a recent paper~\citep{kapidaniArbitraryorderCellMethod2021} by three of the present authors. There a new high order discontinuous Galerkin (DG) method on primal-dual unstructured grids was introduced for the 2D Maxwell equations. The method can be used to efficiently discretise (electromagnetic and acoustic) wave equations in first order form. One of the two unknown fields is discretised in a piecewise conforming way on the original, also called primal, triangulation while the other, which we will refer to as the dual unknown, is piecewise conforming on a barycentric dual, generally polygonal, mesh.
This approach, although unconventional, has several enticing features:
  \begin{enumerate}
    \item There is no need to penalize jumps in the solutions or numerically (i.e., artificially) dissipate energy to achieve spectral correctness of the method (as for some DG approaches, e.g.,~\citep{hesthavenNodalDiscontinuousGalerkin2008}). 
    \item Similarly to DG methods, it leads to  block diagonal mass matrices which provide the amenability to parallelization for HPC implementations.
  \end{enumerate}

On the other hand, open problems remain. The main drawbacks of the previously published work were the following: 
\begin{enumerate}
  \item The blocks in the discrete mass matrix for the dual unknown grow in size with polynomial degree considerably due to the potentially many elements in the original mesh sharing a vertex and the conformity conditions across edges. This is, for example, not the case for the original DG approach which has small fixed size blocks in the mass matrix when increasing the polynomial degree, due to the use of orthogonal polynomials as basis functions.
  \item the local monomial basis used for approximating both primal and dual unknowns is dramatically ill-conditioned for increasing polynomial degree $P$, spoiling $P$--refinement approaches, which is exactly the setting in which DG methods are supposed to shine.
\end{enumerate}

  Fortunately, neither of the aforementioned drawbacks is inherent to the method. Both depend on the choice made for local approximation spaces and their basis functions. By focusing on the quadrilateral mesh generated by the intersection of primal and dual meshes, we can recast the method as a DG-FEM method using quadrilaterals as the basic domain of definition of its approximation spaces. With respect to \citep{kapidaniArbitraryorderCellMethod2021}, we switch from the space $\mathcal{P}^P(K)$ of bivariate polynomials of total degree at most $P$ to the space $\mathcal Q^P(\hat{K})$
  of polynomials of degree at most $P$ separately in each variable on the reference square $\hat{K}$. Subsequently we compose the polynomial basis functions on $\hat K$ with appropriate push-forwards to generic quadrilaterals $K$ to obtain a basis on the physical primal and dual elements.

  After a brief recap of the main ingredients in \cref{sec:lumping}, we show how the new formulations address both issues in \cref{sec:formulations}. 
  Numerical examples validate the improved convergence and efficiency of the new approach in \cref{sec:numerics} and conclusions are thereby drawn in the final section.

\section{Mass lumping on barycentric dual meshes}\label{sec:lumping}
In the present section we provide the notation and basic ingredients for our formulations. We first introduce dual meshes in \cref{sec:meshes}.  Once the quadrilateral nature of the global mesh given by the dual cell method is established, we define polynomial approximation spaces on a reference square in \cref{sec:spaces}. Finally, we provide the recipe to push-forward and glue the local basis functions into global DG spaces similar to the ones presented in \citep{kapidaniArbitraryorderCellMethod2021}.

  \subsection{Barycentric-dual meshes}
  \label{sec:meshes}
  We assume a triangulation $\triangulation$ of the spatial domain $\Omega\subset\R^2$ (in the sense of \cite{Ciarlet}) to be available such that $\cup_{T\in\triangulation} \overline{T}  = \overline{\Omega}$, where the overline denotes the closure of a set. We make the usual assumptions on conforming meshing of discontinuities in the material parameters and call this starting mesh the \emph{primal} mesh when necessity arises.
  \begin{figure}[t!]
    \centering
    \begin{minipage}{0.32\textwidth}
        \centering
        \includegraphics{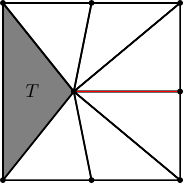}
    \end{minipage}
    \begin{minipage}{0.32\textwidth}
        \centering
        \includegraphics{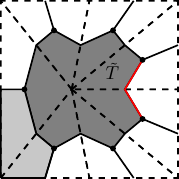}
\end{minipage}
    \begin{minipage}{0.32\textwidth}
        \centering
        \includegraphics{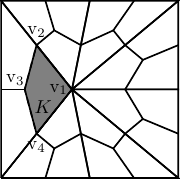}
\end{minipage}
    \caption{From left to right, we first mesh the unit square $\Omega=(0,1)\times(0,1)$ leading to the triangulation $\triangulation$ (we show one triangle $T\in\triangulation$). We then construct the barycentric-dual complex: $\tilde{T}\in\dualmesh$ is dual to a vertex in $\triangulation$. We show in light gray the case of a dual cell corresponding a to a vertex on the boundary of $\Omega$ and in darker gray a dual cell for an internal vertex. We finally highlight the resulting quadrilateral mesh $\quadmesh$ where we emphasize a micro-cell $K$.}
  \label{fig:three_complexes}
  \end{figure}
  We construct in fact a dual mesh for $\triangulation$ which we denote by $\dualmesh$. We perform the construction as in all previous related work on the cell method by taking centroids of triangles as dual vertices and connecting them through segments to the midpoints (centroids) of (primal mesh) edges.
  If a primal edge is part of $\partial\Omega$ there {is} only one such segment originating from the midpoint of the edge. If this is not the case, the union of any two such segments meeting in the centroid of an edge is a poly-line. We call both cases a dual edge, noting that this procedure builds a natural isomorphism between primal and dual edges.
  To complete the construction we obtain a family of two-dimensional sets, called dual cells since they are not necessarily triangles. In fact for each vertex in the interior of $\Omega$ the corresponding dual cell $\tilde{T} \in \dualmesh$ is a (generally non-convex) polygon bounded by dual edges (e.g., the darker gray cell in the second panel of \cref{fig:three_complexes}). For each vertex in $\partial\Omega$ there are instead exactly two edges originating from this vertex. The segments connecting the vertex to their midpoints then complete the boundary of the dual cell (e.g., the lighter gray cell in the second panel of \cref{fig:three_complexes}).
  The barycentric dual procedure outlined above results in a third derived quadrilateral mesh, which is what we are ultimately going to discretise the target partial differential equations with and which we denote by $\quadmesh$. We will also call the quadrilaterals in these mesh micro-cells, to distinguish them from the primal and dual cells.

  We introduce local coordinates for the vertices $\mathbf{v}_K^1$, $\mathbf{v}_K^2$, $\mathbf{v}_K^3$, $\mathbf{v}_K^4$ $\in\mathbb{R}^2$ for each quadrilateral subdomain ${K}\in\quadmesh$ (located e.g., as in the rightmost panel of \cref{fig:three_complexes}) resulting from a non empty intersection of a primal and dual cell. Without loss of generality we always choose $\mathbf v_K^1$ to coincide with a vertex of the primal mesh and the remaining vertices arranged following a counter-clockwise loop around the quadrilateral. This implies that $\mathbf v_K^3$ (which in conjunction with $\mathbf v_K^1$ uniquely identifies the quadrilateral $K$) is always the centroid of a triangle in the primal mesh.
  A physical quadrilateral is then uniquely determined by 
  the continuous and invertible bilinear mapping 
  $\mathbf{F}_K  : \hat K := [0,1]^2 \rightarrow \overline K$,
  that sends vertices of the unit square into vertices of ${K}$, i.e., the unique vector valued mapping $\mathbf{F}_K$
    \begin{align}
      \mathbf F_K(\xi,\eta) :=  (1-\xi)(1-\eta)\mathbf v_K^1+\xi(1-\eta)\mathbf v_K^2+\xi\eta\mathbf v_K^3 + (1-\xi)\eta\mathbf v_K^4,
      \label{eq:F}
    \end{align}
  \noindent which we remark being a standard mapping choice for quadrilateral finite element families (the other being the choice of $[-1,1]\times[-1,1]$ as a reference square domain). 
  We denote by $\mathbf{dF}_K$   the Jacobian matrix associated to $\mathbf{F}_K$, i.e. the $2\times 2$ matrix:
  $$ \mathbf{dF}_K = 
  \begin{pmatrix}
    \\
    (\mathbf v_K^2 - \mathbf v_K^1)(1-\eta) + (\mathbf v_K^3 - \mathbf v_K^4)\eta & 
     (\mathbf v_K^4 - \mathbf v_K^1)(1-\xi)  + (\mathbf v_K^2 - \mathbf v_K^3)\xi\\\,
  \end{pmatrix},
$$
and $J_K = J_K(\xi,\eta)$ its determinant, which is a bilinear polynomial in $\xi$ and $\eta$.

\noindent Since the parametrization $\mathbf{F}_K$ and its inverse are smooth on each quadrilateral, we may define several pushforwards (and pullbacks) for scalar and vector valued functions under change of coordinates from (and to) the unit square to (and from) the physical quadrilateral. 
  The following standard definitions, well known from the FEM literature on de Rham sequences of Finite Element spaces (such as \citep{DemkowiczMonk}) apply locally:
  \begin{equation}
  \begin{array}{ll}
  u(x,y) := \pf_K^{\grad}( \hat u) := \hat u \circ \mathbf{F}_K^{-1} , &\quad \hat u \in \hat{Q}_{P} , \\
  \bm v(x,y) := \pf_K^{\curl}(\hat {\bm v}) := (\mathbf{dF}_K)^{-\top}(\hat {\bm v} \circ \mathbf{F}_K^{-1}) , &\quad \hat{\bm v} \in [\hat{\mathcal Q}_{P}]^2, \\
  \bm w(x,y) := \pf_K^{\div}(\hat {\bm w}) := J_K^{-1} (\mathbf{dF}_K) (\hat{\bm w} \circ \mathbf{F}_K^{-1}) , & \quad \hat{\bm w} \in [\hat{\mathcal Q}_{P}]^2, \label{eq:iotas}
  \end{array}
  \end{equation}
  \noindent where the overhead hat denotes fields defined on the reference square and we denote the space of scalar valued polynomials of degree $P$ with $\hat{\mathcal Q}_{P}$ and $\bm v$, $\bm w$ and $u$ are vector (in boldface) and scalar valued square-integrable functions on $K$.
  The superscripts grad, curl, and div, are due to the fact that the three pushforwards are designed to respectively preserve point values, tangential traces and normal traces of their argument function under changes of coordinates.
  Since the mappings in \cref{eq:iotas} are invertible, we can then easily deduce the expression of the pullbacks through algebraic inversion.

\subsection{Mass-lumped inner products}
\label{sec:inner_products}
In the following we define the approximate (mass lumped) inner products we use in our numerical method.
We start by considering the closed interval $\hat I = [0,1]$ and the Legendre--Gauss--Radau (LGR) quadrature nodes with fixed endpoint in the local variable $\xi\in\hat I$, which consist of $P+1$ points $\{\xi_i\}_{i=0}^P$. They are standard in the literature (e.g. in \citep[Chapter 10.6]{suli_mayers_2003}) and and we assume the points to be sorted in ascending order with $\xi_P=1$. These nodes, with appropriate weights $\{w_i\}_{i=0}^P$, are computed such that integrals of polynomials of degree $2P$ are exactly computable on the unit interval, i.e., if $f(\xi)$ is a polynomial of degree at most $2P$ there holds
\begin{align}
  &\int_0^1 f(\xi) \,\mathrm{d}\xi = \sum_{i=0}^{P} w_{i} f(\xi_i).\label{eq:intrule_1d}
\end{align}

We define a second set of LGR nodes on the interval $[0,1]$, namely the set
$\{\tilde{\xi}_i\}_{i=0}^P$ by $\tilde{\xi}_{P-i} = 1\!-\!\xi_i$.
Through obvious symmetry arguments the second, dual set has the same approximation properties as the primal one for numerical integration, when provided with the corresponding dual weights $\{\tilde{w}_i\}_{i=0}^P$ (in fact $\tilde w_i=w_{P-i}$).

We can use $\{\xi_i\}_{i=0}^P$ to define an interpolatory polynomial basis in a straightforward way: 
the $i$-th Lagrangian polynomial of degree $P$ based on the primal LGR quadrature nodes, denoted as ${\hat{\ell}}_{P}^{(i)}(\xi)$, is defined as the unique polynomial of degree $P$ such that $${\hat{\ell}}_{P}^{(i)}(\xi_j) = \delta_{ij},$$ where $\delta_{ij}$ is the Kronecker delta.
We can likewise define $\tilde {\hat{\ell}}_{P}^{(i)}\in \mathcal P^P(\hat I)$ 
with the analogous Kronecker delta property: $\tilde {\hat{\ell}}_{P}^{(i)}(\tilde \xi_j) = \delta_{ij}$.
Next, we consider the tensorization of these Lagrangian polynomials to extend them from one to two variables (i.e., from $\hat I$ to $\hat K$).
To construct the tensorized Lagrangian polynomials, we define the tensor product of two one-dimensional Lagrangian polynomials ${\hat{\ell}}_{P}^{(i)}(\xi)$ and ${\hat{\ell}}_{P}^{(j)}(\eta)$
\begin{align*}
& {\hat{\ell}}_{P}^{(i,j)}(\xi,\eta) = {\hat{\ell}}_{P}^{(i)}(\xi){\hat{\ell}}_{P}^{(j)}(\eta),\\
& \tilde {\hat{\ell}}_{P}^{(i,j)}(\xi,\eta) = \tilde {\hat{\ell}}_{P}^{(i)}(\xi)\tilde {\hat{\ell}}_{P}^{(j)}(\eta).
\end{align*}
These bivariate polynomials generate the space $$\hat{\mathcal Q}_{P} = \spa\{{\hat{\ell}}_{P}^{(i,j)}(\xi,\eta)\} = \spa\{\tilde {\hat{\ell}}_{P}^{(i,j)}(\xi,\eta)\},$$ which is the space of polynomials of degree up to $P$ in each variable $\xi$, $\eta$.
Similarly one can define vector valued functions:
\begin{align*}
  & \hat {\bm{\ell}}_{P}^{(i,j,k)}(\xi,\eta) = {\hat{\ell}}_{P}^{(i)}(\xi){\hat{\ell}}_{P}^{(j)}(\eta) \hat{\mathbf{e}}_k,\\
  & \tilde {\hat{\bm{\ell}}}_{P}^{(i,j,k)}(\xi,\eta) = \tilde {\hat{\ell}}_{P}^{(i)}(\xi)\tilde {\hat{\ell}}_{P}^{(j)}(\eta)\hat{\mathbf{e}}_k,
\end{align*}
where $\hat{\mathbf{e}}_k$ with $k\in\{1,2\}$ are the unit Euclidean vectors, i.e., $\hat{\mathbf{e}}_1 = (1\; 0)^\top$, $\hat{\mathbf{e}}_2 = (0\; 1)^\top$, and it follows:
\begin{align}
  [\hat{\mathcal Q}_{P}]^2 = \spa\{  {\hat{\bm{\ell}}}_{P}^{(i,j,1)} \} \oplus \spa\{  {\hat{\bm{\ell}}}_{P}^{(i,j,2)} \}=\spa\{  \tilde{{\hat{\bm{\ell}}}}_{P}^{(i,j,1)} \} \oplus \spa\{  \tilde{{\hat{\bm{\ell}}}}_{P}^{(i,j,2)} \},
   \label{eq:QPhat}
\end{align}
where $\oplus$ denotes the direct sum of vector spaces.
We remark that, extending the integration rule of \eqref{eq:intrule_1d} to two dimensions we can define the approximate inner product, for both scalar functions $\hat f, \hat g \in C(\hat{K})$ and vector valued functions $\hat{\bm u}, \hat{\bm v} \in [C(\hat{K})]^2$, given by 
\begin{align}
  & \left< \hat f,\hat g \right>_{\hat{K}}^{P} := \sum_{i=0}^{P} \sum_{j=0}^{P}w_{i} w_{j} \hat f(\xi_i,\xi_j) \hat g(\xi_i,\xi_j) \approx \int_0^1\!\!\mathrm{d}\eta\int_0^1\!\!\mathrm{d}\xi\,\hat f(\xi,\eta) \hat g(\xi,\eta)  ,\label{eq:intrule_2d_scalar} \\
  & \left< \hat{\bm u}, \hat{\bm v} \right>_{\hat{K}}^{P} := \sum_{i=0}^{P} \sum_{j=0}^{P}w_{i} w_{j} \hat {\bm u} (\xi_i,\xi_j) \cdot \hat {\bm v} (\xi_i,\xi_j) \approx \int_0^1\!\!\mathrm{d}\eta\int_0^1\!\!\mathrm{d}\xi\,  \hat {\bm u} (\xi,\eta) \cdot \hat {\bm v} (\xi,\eta){,} \label{eq:intrule_2d_vec}
\end{align}
\noindent where the same notation will be used for inner products computed using dual integration rules. This is a slight abuse of notation, but the choice of the integration rule based on $\{\tilde{\xi}_i\}_{i=0}^P$ rather than $\{{\xi}_i\}_{i=0}^P$ will be obvious from the context in \cref{sec:formulations}. We call these inner products \emph{mass-lumped} inner products. 
Likewise the following properties which we present for the primal integration nodes can be stated nearly verbatim for the dual ones. Therefore in the following paragraphs we only go into detail for the primal integration rules.

On the reference square it is easy to see that the following exact orthogonality properties hold:
\begin{align*}
  \left< {\hat{\ell}}_{P}^{(i,j)},{\hat{\ell}}_{P}^{(l,m)} \right>_{\hat{K}}^{P} =\; & w_i w_j\delta_{il}\delta_{jm}, \\
  \left< \hat {\bm{\ell}}_{P}^{(i,j,k)},\hat {\bm{\ell}}_{P}^{(l,m,n)} \right>_{\hat{K}}^{P} =\; & w_i w_j\delta_{il}\delta_{jm}\delta_{kn},
\end{align*}

\noindent where $\delta$ is again the Kronecker symbol. We remark that the inner products above are actually identical to the $\mathrm{L}^2$ inner products, since the integrands are polynomials of degree $2P$. Both properties will be pivotal in the following.

We can thus define how our approximate inner products then look like on the physical quadrilateral $K$ in the mesh for scalar continuous functions $f,g$ on $K$, based on LGR integration rules 
\begin{align}
  \left< f,g\right>_K^P:=\left< J_K f\circ \mathbf F_K , g\circ \mathbf F_K\right>^P_{\hat K},
  \label{eq:defn_lumped}
\end{align}
with the obvious same definition for vector-valued $\left<\bm u,\bm v\right>_{\hat K}^P$ and continuous $\bm u,\bm v$.

In the particular case in which $u,\bm v,\bm w$ (and $u',\bm v',\bm w'$) are push-forwards as in \cref{eq:iotas} we obtain
\begin{align*}
  \left< u, u'\right>_{K}^{P} & =\left<  \pf_K^{\grad}( \hat u)\circ\mathbf F_K,\pf_K^{\grad}( \hat u')\circ\mathbf F_K \right>_{\hat{K}}^{P}
  =\left< \hat u,J_K\hat u'\right>_{\hat{K}}^{P},\\
  \left< \bm v,\bm v' \right>_{K}^{P} & =\left< \pf_K^{\curl}( \hat {\bm v})\circ\mathbf F_K,\pf_K^{\curl}({\hat {\bm v}}')\circ\mathbf F_K \right>_{\hat{K}}^{P}
  =\left< \hat {\bm v},\mathbb{G}^K{\hat {\bm v}}'\right>_{\hat{K}}^{P},\\
  \left< \bm w,\bm w' \right>_{K}^{P} & =\left< \pf_K^{\div}( \hat {\bm w})\circ\mathbf F_K,\pf_K^{\div}({\hat {\bm w}}')\circ\mathbf F_K \right>_{\hat{K}}^{P}
  =\left< \hat {\bm w},\mathbb{H}^K{\hat{\bm w}}'\right>_{\hat{K}}^{P},
\end{align*}
\noindent where we have absorbed the (smooth) Jacobian determinant of the mapping within one inner product argument without loss of generality.
We furthermore used the symmetric matrices:
\begin{align*}
  &\mathbb{G}^K = (\mathbf{dF}_K)^{-1} J_K (\mathbf{dF}_K)^{-\top},\\
  &\mathbb{H}^K = (\mathbf{dF}_K)^{\top} J_K^{-1} (\mathbf{dF}_K),
\end{align*}
where we have omitted for readability the dependence of all quantites on $\xi$ and $\eta$ and it is clear that $\mathbb{G}^K \mathbb{H}^K = \mathbb{I}^{2\times2}$.
The following (quasi) orthogonality conditions then hold:
\begin{align*}
  \left< \pf_K^{\grad}({\hat{\ell}}_{P}^{(i,j)}),\pf_K^{\grad}({\hat{\ell}}_{P}^{(l,m)}) \right>_{K}^P =\; & w_i w_j\delta_{il}\delta_{jm} J_K (\xi_i,\xi_j), \\
  \left< \pf_K^{\curl}(\hat {\bm{\ell}}_{P}^{(i,j,k)}),\pf_K^{\curl}(\hat {\bm{\ell}}_{P}^{(l,m,n)}) \right>_{K}^{P} =\; & w_i w_j\delta_{il}\delta_{jm} (\mathbb{G}^K)_{kn}(\xi_i,\xi_j),\\
  \left< \pf_K^{\div}(\hat {\bm{\ell}}_{P}^{(i,j,k)}),\pf_K^{\div}(\hat {\bm{\ell}}_{P}^{(l,m,n)}) \right>_{K}^{P} =\; & w_i w_j\delta_{il}\delta_{jm} (\mathbb{H}^K)_{kn}(\xi_i,\xi_j), 
\end{align*}
where $(\mathbb{G}^K)_{kn}(\xi_i,\xi_j)$ is the $(k,n)$-th entry of the matrix $\mathbb{G}^K$ evaluated at coordinates $(\xi_i,\xi_j)$ and $(\mathbb{H}^K)_{kn}(\xi_i,\xi_j)$ is the $(k,n)$-th entry of the matrix $\mathbb{H}^K$ evaluated at coordinates $(\xi_i,\xi_j)$. Since $k,n\in\{1,2\}$, assembling the last two above expressions into inner-product matrices yields 2-by-2 block diagonal matrices.



\subsection{Nodal spaces on dual cells}
\label{sec:spaces}
Now that we have defined the mass-lumped inner products we focus on the construction of the (local) bases of our approximation spaces. Again the construction is very similar for the primal and dual bases. Since it is the more exotic case, we focus on the construction of the dual spaces.

We choose the degrees of freedom of our spaces to be point values of the respective fields in the integration points defined above to exploit the orthogonalites derived in the previous subsection.
To this end we take a dual cell ${\tilde T} \in \dualmesh$. This is the union of a finite number of quadrilaterals $K$ on which we construct three types of discrete spaces $\Xgradd{\tilde T}$, $\Xcurld{\tilde T}$, $\Xdivd{\tilde T}$ spanned by piecewise polynomials.
\subsubsection{Basis functions for $\Xgradd {\tilde T}$}
\begin{figure}[th]
  \begin{subfigure}{0.5\textwidth}
    \centering
    \includegraphics{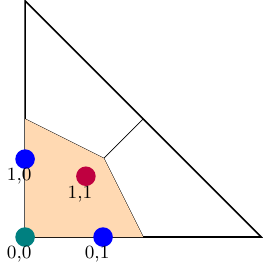}
    \caption{$P=1$}
    \label{subfig:basis_h1dual_1_K}
  \end{subfigure}
  \begin{subfigure}{0.5\textwidth}
    \centering
    \includegraphics{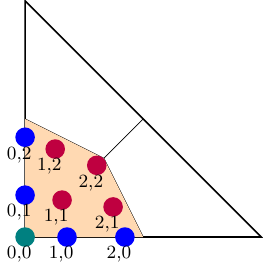}
    \caption{$P=2$}
    \label{subfig:basis_h1dual_2_K}
  \end{subfigure}
  \caption{Points associated to the degrees of freedom (DoFs) of $X_1^{\grad}(\tilde{\mathcal T_h})$ and $X_2^{\grad}(\tilde{\mathcal T_h})$ on one quadrilateral. Red dots correspond to basis functions $u^K_i$, blue ones to $u_i^E$ and teal ones to $u^V$.}
  \label{fig:basis_h1dual_K}
\end{figure}

\begin{figure}[th]
    \centering
  \begin{subfigure}{0.4\textwidth}
    \centering
    \includegraphics[width=\textwidth]{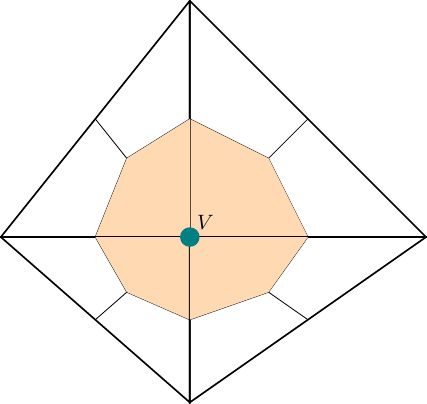}
    \caption{location of the vertex DoF}
    \label{subfig:basis_h1dual_vertex}
  \end{subfigure}
  \begin{subfigure}{0.4\textwidth}
    \centering
    \includegraphics[width=\textwidth]{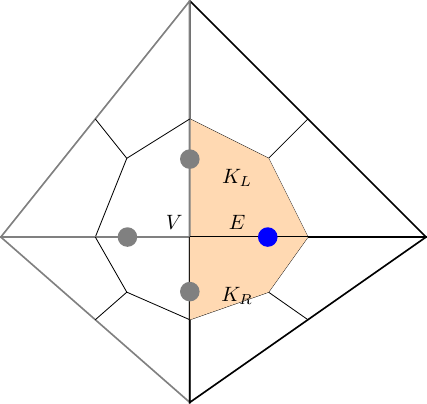}
    \caption{location of the edge DoFs}
    \label{subfig:basis_h1dual_edge}
  \end{subfigure}\\

  \begin{subfigure}{0.4\textwidth}
    \centering
    \includegraphics[width=\textwidth]{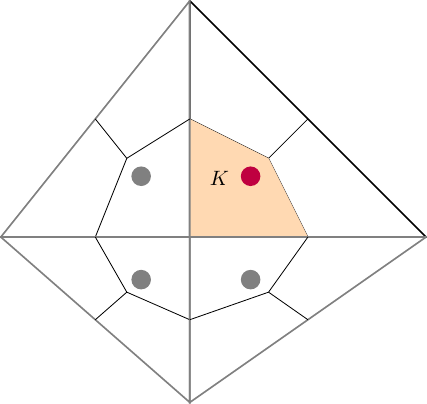}
    \caption{location of the face DoFs}
    \label{subfig:basis_h1dual_face}
  \end{subfigure}
  \begin{subfigure}{0.4\textwidth}
    \centering
    \includegraphics[width=\textwidth]{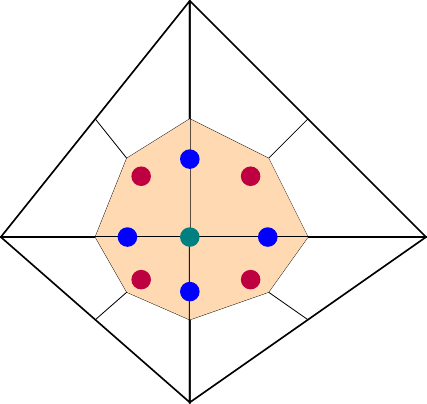}
    \caption{location of all DoFs}
    \label{subfig:basis_h1dual_all}
  \end{subfigure}
  \caption{Points associated to the vertex, edge and face DoFs \cref{subfig:basis_h1dual_vertex,subfig:basis_h1dual_edge,subfig:basis_h1dual_face} of $X_1^{\grad}(\tilde{\mathcal T_h})$ on one element of the dual mesh (\cref{subfig:basis_h1dual_all} colored). The colored quadrilaterals mark the support of the basis function corresponding to the colored DoF, where the same color coding as in \cref{fig:basis_h1dual_K} applies.}
  \label{fig:basis_h1dual}
\end{figure}
We start by remarking that for a given dual cell ${\tilde T}$ there exists exactly one vertex $V$ of the primal such that $V\in\overline{{\tilde T}}$. Then the function which is interpolatory at $V$, defined by
$$u^{V} (\bm x)= 
\begin{cases}
  \pf_K^{\grad}({\hat{\ell}}_{P}^{(0,0)}) & \text{if } \bm x \in K \subset {\tilde T} \quad \forall K \text{ s.t. } V\in\partial K, \\
  0 & \text{otherwise,} 
\end{cases}
$$
is continuous on ${\tilde T}$ (cf., \cref{fig:basis_h1dual_K}, DoFs with index $(0,0)$ and \cref{subfig:basis_h1dual_vertex}).
Additionally, we have edge functions: for each edge $E$ originating from $V$ in the skeleton of $\quadmesh$, there are two quadrilaterals $K_L$, $K_R$ for which $E\subseteq \partial K_L$ and $E\subseteq\partial K_R$ (cf., Figure \cref{subfig:basis_h1dual_edge}), where the subscripts stand for left and right respectively and are motivated by the right-handed corkscrew rule providing an inner orientation for ${\tilde T}$. There are then $P$ basis functions of the kind:
  $$u_i^{E} (\bm x)= 
  \begin{cases}
    \pf_{K_L}^{\grad}({\hat{\ell}}_{P}^{(0,i)}) & \text{if } \bm x \in K_L, \\
    \pf_{K_R}^{\grad}({\hat{\ell}}_{P}^{(i,0)}) & \text{if } \bm x \in K_R, \\
    0 & \text{otherwise,} 
  \end{cases}
  $$
  with $i=1,\dots,P$. We remark that the definition trivially extends to the case in which $E\subset\partial\Omega$ where we simply have either $K_L = \emptyset$ or $K_R = \emptyset$.
  
  Finally, for each quadrilateral $K\subset{\tilde T}$ there are $P^2$ functions of the kind:
  $$u_{(j-1)P + i}^{K} (\bm x)= 
  \begin{cases}
    \pf_K^{\grad}({\hat{\ell}}_{P}^{(i,j)}) & \text{if } \bm x \in K \subset {\tilde T}, \\
    0 & \text{otherwise.} 
  \end{cases}
  $$
  with $i,j=1,\dots,P$. We remark that these basis functions are compactly supported on a single quadrilateral $K$ (cf., \cref{subfig:basis_h1dual_face}).

We then define $\Xgradd {\tilde T}$ as the span of the union of the three sets of basis functions above and remark that $\Xgradd {\tilde T}\subset \mathrm{H}^1({\tilde T})$ where $\mathrm{H}^1({\tilde T})$ is the space of square integrable functions on ${\tilde T}$ with square integrable gradient, i.e.,
\begin{align*} 
  \mathrm{H}^1({\tilde T}) = \left\{ u \in \mathrm{L}^2({\tilde T}) : \grad u \in \mathrm{L}^2({\tilde T}) \right\}. 
\end{align*}

A pictorial representation of the degrees of freedom {on a single quadrilateral} for $P=1$ and $P=2$ is given in \cref{fig:basis_h1dual_K} {and for $P=1$ on a single dual element in \cref{fig:basis_h1dual}}.



\subsubsection{Basis functions for $\Xcurld {\tilde T}$}
For the vector valued space $\Xcurld {\tilde T}$ we have the following classification of basis functions:
\begin{itemize}
  \item For each edge $E\subseteq {\tilde T} \setminus {\partial {\tilde T}}$, again there are two quadrilaterals $K_L$, $K_R$ for which $E\subseteq \partial K_L$ and $E\subseteq\partial K_R$, as above. There are then $P+1$ basis functions of the kind:
    $$\bm u_{i+1}^{E} (\bm x)= 
  \begin{cases}
    \pf_{K_L}^{\curl}({\hat{\bm \ell}}_{P}^{(0,i,2)}) & \text{if } \bm x \in K_L, \\
    \pf_{K_R}^{\curl}({\hat{\bm \ell}}_{P}^{(i,0,1)}) & \text{if } \bm x \in K_R, \\
    0 & \text{otherwise,} 
  \end{cases}
  $$
  with $i=0,\dots,P$ (note we start from $0$). Again the definition trivially extends to the case in which $E\subset\partial\Omega$ where we have simply either $K_L = \emptyset$ or $K_R = \emptyset$.
  \item For each quadrilateral $K\subset{\tilde T}$ there are $P(P+1)$ functions of the kind:
    $$\bm u_{(j-1)(P+1) + i+1}^{K} (\bm x)= 
  \begin{cases}
    \pf_K^{\curl}({\hat{\bm \ell}}_{P}^{(i,j,1)}) & \text{if } \bm x \in K \subset {\tilde T}, \\
    0 & \text{otherwise,} 
  \end{cases}
  $$
  with $i=0,\dots,P$ and $j=1,\dots,P$ and additional $P(P+1)$ functions of the kind:
    $$\bm u_{(j-1+P)(P+1) + i+1}^{K} (\bm x)= 
  \begin{cases}
    \pf_K^{\curl}({\hat{\bm \ell}}_{P}^{(i,j,2)}) & \text{if } \bm x \in K \subset {\tilde T}, \\
    0 & \text{otherwise,} 
  \end{cases}
    $$
    with $i=1,\dots,P$ and $j=0,\dots,P$. Both kind{s} of functions above are compactly supported on a single quadrilateral $K$.
\end{itemize}

\begin{figure}
    \centering
  \begin{subfigure}{0.4\textwidth}
    \includegraphics[width=\textwidth]{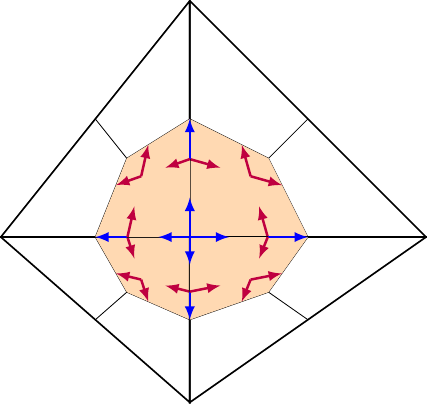}
    \centering
    \caption{$P=1$}
    \label{subfig:basis_hcurldual_1}
  \end{subfigure}
  \begin{subfigure}{0.4\textwidth}
    \includegraphics[width=\textwidth]{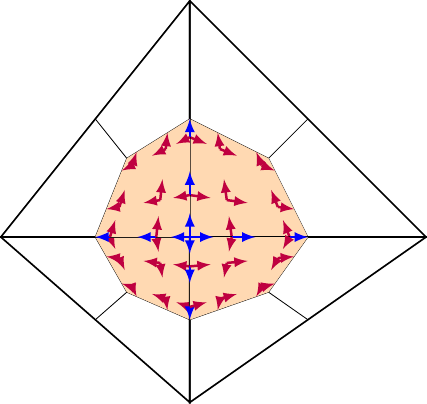}
    \centering
    \caption{$P=2$}
    \label{subfig:basis_hcurldual_2}
  \end{subfigure}
  \caption{Distribution of the DoFs of $\tilde X_1^{\curl}(\dualmesh)$ and $\tilde X_2^{\curl}(\dualmesh)$ on one dual element (colored).}
  \label{fig:basis_hcurldual}
\end{figure}

We then define $\Xcurld{\tilde{T}}$ as the span of the union of the three sets of basis functions above. A pictorial representation of the degrees of freedom for $P=1$ and $P=2$ is given in \cref{fig:basis_hcurldual}.
We again remark that $\Xcurld{{\tilde T}}$ is a subspace of the space $H(\curl;{\tilde T})$ defined as follows:
\begin{align*} 
  & \mathrm{H}({\curl};{\tilde T}) = \left\{ \bm{v} \in \mathrm{L}^2({\tilde T})^2 : \curl \bm{v} \in \mathrm{L}^2({\tilde T})^2 \right\}.
\end{align*}
{T}he basis functions for $\Xdivd{\tilde T}$ are constructed analogously to the ones for $\Xcurld{\tilde T}$ and we omit the details here. It suffices to swap Cartesian components in edge based functions and swap $i$, $j$ iterating indices in basis functions which are instead locally supported on one quadrilateral $K$.
We remark that $\Xdivd{\tilde T}$ is a subspace of the space $H(\div;{\tilde T})$ defined as follows:
\begin{align*} 
  H(\div;{\tilde T}) &= \left\{ \bm{u} \in \mathrm{L}^2({\tilde T})^2 : \div \bm{u} \in \mathrm{L}^2({\tilde T}) \right\}.
\end{align*}

\subsubsection{Construction of the global spaces}
\begin{figure}[th!]
 \begin{subfigure}{0.5\textwidth}
  \centering
   \includegraphics{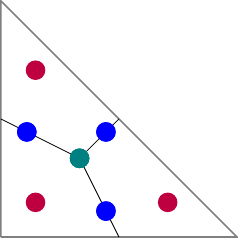}
   \caption{{$P=1$}}
   \label{subfig:basis_h1primal_1}
 \end{subfigure}
 \begin{subfigure}{0.5\textwidth}
  \centering
   \includegraphics{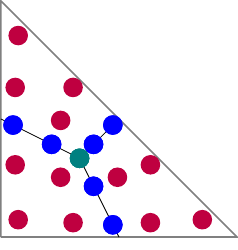}
   \caption{{$P=2$}}
   \label{subfig:basis_h1primal_2}
 \end{subfigure}
  \caption{Points associated to the DoFs of $\Xgradp\triangulation$ for $P=1,2$ on one primal element.}
 \label{fig:basis_h1primal}
\end{figure}

\begin{figure}[th!]
 \begin{subfigure}{0.5\textwidth}
  \centering
   \includegraphics{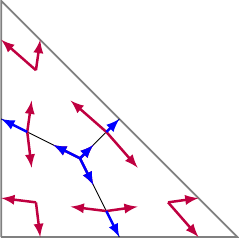}
   \caption{{$P=1$}}
   \label{subfig:basis_hcurlprimal_1}
 \end{subfigure}
 \begin{subfigure}{0.5\textwidth}
  \centering
   \includegraphics{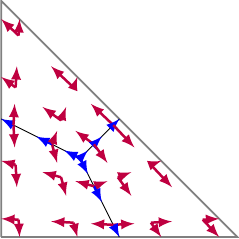}
   \caption{{$P=2$}}
   \label{subfig:basis_hcurlprimal_2}
 \end{subfigure}
  \caption{Points associated to the DoFs of $\Xcurlp\triangulation$ for $P=1,2$ on one primal element.}
 \label{fig:basis_hcurlprimal}
\end{figure}

Finally, to obtain the discrete spaces on the whole mesh, we define the global space $\Xgradd \dualmesh$ as Cartesian product of the local spaces $\Xgradd {\tilde T}$ for all ${\tilde T} \in \dualmesh$, and similarly for the vector valued fields:
\begin{equation}
  \begin{array}{ll}
    \Xgradd \dualmesh := & \prod_{{\tilde T} \in \dualmesh} \Xgradd {\tilde T}, \\
    \Xcurld \dualmesh := & \prod_{{\tilde T} \in \dualmesh} \Xcurld{\tilde T},\\
    \Xdivd \dualmesh  := & \prod_{{\tilde T} \in \dualmesh} \Xdivd {\tilde T},
  \end{array}
\end{equation}
and degrees of freedom of the global spaces are obtained by the union of the degrees of freedom of the local spaces. We omit the explicit construction of bases for the spaces $\Xgradp T$,  $\Xcurlp T$, $\Xdivp T$ in the case of the primal triangles $T\in\triangulation$ since their construction is analogous to what is done above for the case of the dual mesh. We still provide the pictorial representation of the degrees of freedom of the scalar valued space in \cref{fig:basis_h1primal} and the vector valued one in \cref{fig:basis_hcurlprimal} for a single triangle in the mesh.


\section{Discrete formulations for waves in 2D}\label{sec:formulations}
We will now use the notation and tools from the \cref{sec:lumping} to show how to efficiently solve the acoustic and electromagnetic wave equations via the dual cell method in 2D. We will devote specific attention to the similarities and the differences between the mass lumped approach here and the already published approach in \citep{kapidaniArbitraryorderCellMethod2021}.

\subsection{Maxwell equations} \label{sec:maxwell_eq}
Consider the time-dependent Maxwell equations in their first-order form for an electric field $\bm{E}(\bm x, t)\in \mathrm{L}^2(\mathrm{H}({\curl};\Omega);[0,T])$ and a magnetic field $\bm{H}(\bm x,t) \in \mathrm{L}^2(\mathrm{H}({\rot};\Omega);[0,T])$, such that:
\begin{align}
    &\varepsilon\frac{\partial \bm{E}}{\partial t} = \rot {H} - \bm{J}, & \text{in } \Omega \times (0, T), \label{eq:maxwell1} \\
    &\mu\frac{\partial {H}}{\partial t} = - \curl \bm{E},  & \text{in } \Omega \times (0, T), \label{eq:maxwell2} \\
    &\bm{E} \times \normal_\Omega = 0, & \text{on } \partial\Omega \times (0, T), \label{eq:maxwellbc} \\
    &\bm{E}(\bm x, 0) = \bm{E}_0(\bm x)\in \mathrm{H}({\curl};\Omega), & H(\bm x, 0) = {H}_0(\bm x)\in \mathrm{H}({\rot};\Omega) & \text{ in } \Omega, \label{eq:maxwelinit}
\end{align}
where $\bm{E}_0(\bm x)\in \mathrm{H}({\curl};\Omega)$, ${H}_0(\bm x)\in \mathrm{H}({\rot};\Omega)$, and where $\varepsilon,\mu \in L^\infty(\Omega)$ are the dielectric permittivity and magnetic permeability, assumed to be scalar valued and time-invariant for simplicity of exposition. The vector field $\bm{J}(\bm x, t)$ (assumed to be $0$ in most of our numerical experiments) is the electric current density (a suitable a priori known source of the sytem) and $T>0$ is the final simulation time.
This is the system on which some of the authors also focused in \citep{kapidaniArbitraryorderCellMethod2021} and provides also the most challenging test case for the new approach, since mass lumping for high order versions of edge elements is not as straightforward to achieve as for the scalar valued pressure in the acoustic case.
We point out that \cref{eq:maxwell1}--\cref{eq:maxwelinit} is obtained using the transverse magnetic ansatz for the actual vector valued magnetic field (necessarily defined in a three-dimensional domain) to be of the form $\bm{H}=H\hat{\bm z}$, where $\hat{\bm z}$ is the unit vector in the $z$-direction and $H\in H(\rot;\Omega)$.
By Galerkin testing and integration by parts of \eqref{eq:maxwell1}--\eqref{eq:maxwell2} (and neglecting current sources for the sake of brevity) on each quadrilateral~$K\in\quadmesh$, boundary integrals appear on the segments bounding $K$. We obtain the semi-discrete weak formulation to find $\bm{E}_h\in{\tilde{X}}^{{\curl}}_{P}(\dualmesh)$ and $H_h \in X_{P}^{\rot}(\triangulation)$ such that:
\begin{align}
   \sumoveralldtrigs\sumoverallelementsindtrig \left< \varepsilon \frac{\partial \bm{E}_h}{\partial t}, \bm e  \right>_{K}^{P}\!=\!
    \sumoveralldtrigs\left(\sumoverallelementsindtrig  \int_{K} H_h\hat{\bm z}\cdot\curl \bm e + \sumoverallfacetslocald \int_{F} \bm e \cdot H_h \hat{\bm z} \times \normal_K \right), \label{eq:uw_am} \\
    \sumoverallptrigs\sumoverallelementsinptrig \left< \mu \frac{\partial H_h}{\partial t},h \right>_{K}^{P}\!=\!
    \sumoverallptrigs\left(\sumoverallelementsinptrig -\int_{K} \bm{E}_h\cdot\rot {h} + 
    \sumoverallfacetslocalp \int_{F} \bm{E}_h \cdot h\hat{\bm z} \times \normal_K \right), \label{eq:uw_faraday}
\end{align}
\noindent holding for all $\bm e \in {\tilde{X}}^{{\curl}}_{P}(\dualmesh)$, and $h \in X_{P}^{\rot}(\triangulation)$, where $\normal_K$ denotes the outer normal on each element boundary. The above system implies weak imposition of the $\bm{E}_h\times\normal_\Omega = 0$ boundary conditions. We also require $\bm{E}_h|_{t=0}$ and $H_h|_{t=0}$ to fulfill the initial conditions through an $\mathrm{L}^2$-projection at the initial time.
We close the subsection with a final important statement regarding the Maxwell system, which clarifies that the most succint way of writing its weak formulation, chosen in the manuscript, is not necessarily the most efficient one for the practical implementation.
\begin{remark}
  \label{rmk:impl}
  For the implementation of the presented method one may apply integration by parts to \cref{eq:uw_am}. The reader may easily verify that in this case all the edge contributions which are not part of $\partial {\tilde T}$ for any ${\tilde T}\in\tilde{\mathcal T}$ cancel out. Thus, one is only left with integral contributions on either a triangle $T$, or its boundary $\partial T$, i.e., on the original triangulation assumed available from a FEM perspective. An implementation can therefore exploit existing finite element codes (as we did within NetGen/NGSolve, \citep{netgen,ngsolve}),
without explicitly having to generate the dual mesh.
\end{remark}

\subsection{Acoustic wave equations}
In this section, we briefly sketch how the same method can be applied to the simulation of acoustic waves, without any practical differences in computational efficiency. 
Analogously to \cref{sec:maxwell_eq} we start from a strong formulation in the velocity-pressure first-order system form, i.e.{,} the initial boundary value problem of finding $\bm{V}(\bm x, t)\in \mathrm{L}^2(\mathrm{H}({\div};\Omega);[0,T])$ and $Q(\bm x,t) \in \mathrm{L}^2(\mathrm{H}({\grad};\Omega);[0,T])$ such that:
\begin{align}
    &\frac{\partial Q}{\partial t} = \rho_0 c^2 \div \bm{V}+f, & \text{in } \Omega \times (0, T), \label{eq:acoustic1} \\
    &\rho_0\frac{\partial \bm{V}}{\partial t} = \grad Q, & \text{in } \Omega \times (0, T), \label{eq:acoustic2} \\
    &\bm{V} \times \normal = 0 & \text{on } \partial\Omega \times (0, T), \label{eq:acousticbc} \\
    &\bm{V}(\bm r, 0) = \bm{V}_0(\bm r), & \bm{Q}(\bm r, 0) = \bm{Q}_0(\bm r) & \text{ in } \Omega, \label{eq:acousticinit}
\end{align}
with $\bm{V}\cdot\normal_\Omega=0$, where \(Q\) is the acoustic pressure, \(\bm{V}\) is the particle velocity, \(\rho_0\) is the reference density, \(c\) is the speed of sound and $f(\bm x ,t)$ is a suitable given forcing pressure. We again have initial conditions $\bm{V}_0(\bm x)\in \mathrm{H}({\div};\Omega)$, ${Q}_0(\bm x)\in \mathrm{H}({\grad};\Omega)$.
The semi-discrete weak formulation of \cref{eq:acoustic1}--\cref{eq:acoustic2} then seeks \(Q_h \in X_P^{\grad}(\triangulation)\) and \(\bm{V}_h \in {\tilde X}_P^{\div}(\dualmesh)\) such that for all test functions \(Q \in X_P^{\grad}(\triangulation)\) and \(\bm{V} \in{\tilde X}_P^{\div}(\dualmesh)\):

%
\begin{align*}
     \sumoveralldtrigs\sumoverallelementsindtrig   \left< \frac{1}{\rho_0 c^2} \frac{\partial Q_h}{\partial t}, q \right>_{K}^{P}\!&=&\!
         \sumoveralldtrigs \left(\sumoverallelementsindtrig \int_{K} \bm{V}_h \cdot \grad q + \sumoverallfacetslocald \int_{F} q\; (\bm{V}_h \cdot \normal_K) \right), \\
         \sumoverallptrigs\sumoverallelementsinptrig\left< \rho_0 \frac{\partial \bm{V}_h}{\partial t}, \bm{v} \right>_{K}^{P}\!&=&\!
         \sumoverallptrigs\left(\sumoverallelementsinptrig  \int_{K} Q_h \; \div \bm{v} + 
         \sumoverallfacetslocalp \int_{F} Q_h\; (\bm{v} \cdot \normal_K) \right), 
\end{align*}
holding for all $\bm v \in {\tilde{X}}^{{\div}}_{P}(\dualmesh)$, and $q \in X_{P}^{\grad}(\triangulation)$, where $\normal_K$ denotes again the outer normal on each boundary of $K$. 
The above system implies weak imposition of the $\bm{V}_h\cdot\normal_\Omega = 0$ boundary conditions. We again also require $\bm{V}_h|_{t=0}$ and $Q_h|_{t=0}$ to fulfill the initial conditions through an $\mathrm{L}^2$-projection at the initial time.
Finally we stress that \cref{rmk:impl} (with appropriate modifications) also hold for the acoustic case.

\subsection{Related methods}
Before putting the new mass lumped approach to the numerical test we dedicate a subsection to address the relationship, similarities and differences between the presented method and related approaches available in the literature.
The proposed approach is, as was the case for \citep{kapidaniArbitraryorderCellMethod2021}, based on the seminal work on low order barycentric dual grid methods in \citep{codecasaExplicitConsistentConditionally2008,codecasaNovelFDTDTechnique2018,cicuttinGPUAcceleratedTimeDomain2018,kapidaniTimeDomainCellMethod2020}. The barycentric dual is not the only choice that provides a possible extension to high order. In fact Chung and co-authors have championed staggered DG methods in several works on structured \citep{chungConvergenceSuperconvergenceStaggered2013} and unstructured \citep{chungOptimalDiscontinuousGalerkin2009,b.cockburnDiscontinuousGalerkinMethods2012,chungStaggeredDGMethod2014} grids, to obtain a penalty free approach which presents block-diagonal mass matrices. 
Nevertheless, in the very relevant unstructured case their underlying micro-cells, on which one exploits the same integration by parts formula are still triangles, since the Worsey--Farin split mesh is used. This can be a drawback in two ways: firstly, since each triangle corner is bisected, the shape-regularity constants degrade with respect to the original triangular mesh worsening the system's conditioning. In our case the shape-regularity constants of the new mesh are the same as the ones given by the initial FE mesh and less stringent than the ones of the Worsey--Farin split \citep{gongNoteShapeRegularity2023}.
Secondly, improving the orthogonality of basis functions through mass lumping techniques amounts to finding an integration rule which allows an $\mathrm{H}({\curl},T)$-conforming (or $\mathrm{H}({\div},T)$-conforming for mixed formulation of the acoustic case) nodal basis on a Worsey--Farin split of a triangle to be constructed, which is ultimately just as hard as finding such a construction for a fully $\mathrm{H}({\curl},\Omega)$-conforming discrete space on the whole triangulation. Even though progress has been made recently in this regard with mass lumping schemes for acoustic and electromagnetic wave equations of higher order (see \citep{geeversNewHigherOrderMassLumped2018,eggerMasslumpedMixedFinite2020,eggerSecondOrderFiniteElement2021}), a general recipe for arbitrary polynomial degrees remained an open question. Our mass-lumping technique is formulated on a quadrilateral mesh and it therefore bypasses the issue completely by relying on tensor-product integration rules.


\subsection{Time Discretisation}
\label{sec:time_disc}
We turn our focus back {to} the Maxwell equations as an example for the time discretisation, since the wave equation is treated in a completely analogous way.
The common leap frog scheme for the time discretisation of \cref{eq:uw_am}, \cref{eq:uw_faraday} is used as a convenient choice to test the new mass lumping approach, since it is a very easy to implement symplectic integrator and it is second order accurate in time. If we denote with ${\mathbf{h}}\in\mathbb{R}^{\dim(X_{P}^{\rot}(\triangulation))}$ and $\mathbf{e}\in\mathbb{R}^{\dim({\tilde{X}}^{{\curl}}_{P}(\dualmesh))}$ the vectors of degrees of freedom for the magnetic and electric field respectively, appropriately ordered to evidence the block diagonal structures of the mass matrices, the scheme is given by the following update rules:

  \begin{align*}
  & {\mathbf{h}}^{1/2} = { \mathbf{h}}^{0} -\frac{\Delta t}{2} \,\mathbf{M}_{{\mu}}^{-1} \mathbf{C} \mathbf{e}^0,\\
  & \mathbf{e}^{1} = \mathbf{e}^{0} +\Delta t \,\mathbf{M}_{{\varepsilon}}^{-1} \mathbf{C}^\top {\mathbf{h}}^{1/2},
  \end{align*}
for the first update of the magnetic field and 
\begin{align*}
& {\mathbf{h}}^{n+1/2} = {\mathbf{h}}^{n-1/2} -\Delta t \,\mathbf{M}_{{\mu}}^{-1} \mathbf{C} \mathbf{e}^n,\\
& \mathbf{e}^{n+1} = \mathbf{e}^{n} +\Delta t \,\mathbf{M}_{{\varepsilon}}^{-1} \mathbf{C}^\top {\mathbf{h}}^{n+1/2},
\end{align*}
for all $n\geq1$, in which $\mathbf{e}^n$ is the approximation computed at time instant $n\,\Delta t$ and ${\mathbf{h}}^{n+1/2}$ is computed at time instant $(n+1/2)\,\Delta t$, for $n=0,\dots,N$. The final time $T$ is divided into $N$ time steps of size $\Delta t = T/N$. Sparse matrix $\mathbf{C}$ is the matrix representation of the discrete curl operator, which is the r.h.s. of \cref{eq:uw_faraday}, while $\mathbf{M}_{{\mu}}$ and $\mathbf{M}_{{\varepsilon}}$ are the mass matrices for the magnetic and electric field respectively.

What we want to stress about the fully discrete formulation is that the l.h.s. gives rise to block diagonal system matrices matrices $\mathbf{M}_{{\mu}}$ and $\mathbf{M}_{{\varepsilon}}$. They thus have block diagonal
inverses where the block size does not increase with increasing polynomial degree, differently from what was achieved in previous work were scaled monomials were used as a basis.

In order to guarantee stability, the time step $\Delta t$ is computed in such a way that 
\begin{align}
  \Delta t < t_0:=2 / \sqrt{\lambda_M},
  \label{eq:cfl}
\end{align}
where $\lambda_M$ is the maximum of the eigenvalues $\lambda$ of the eigenvalue problem
\begin{align}
& \mathbf{C} \mathbf{M}_{{\varepsilon}}^{-1} \mathbf{C}^\top {\mathbf{h}} = \lambda \mathbf{M}_{{\mu}} {\mathbf{h}},
  \label{eq:disc_evp}
\end{align} 
for the eigenvectors ${\mathbf{h}}$. 
The dependence of this so called Courant-Friedrichs-Lewy (CFL) condition on the discretization parameters is studied numerically in \cref{sec:cfl}.

\section{Numerics}\label{sec:numerics}
We perform several numerical experiments to underline the applicability, efficiency and convergence of our method. {To this end we use an implementation of our method in the high-order finite element library Netgen/NGSolve \cite{netgen,ngsolve} along the lines of \cref{rmk:impl}}.

In particular, in \cref{sec:disc_evp} we show that we obtain spectral convergence of the underlying eigenvalue problem with polynomial convergence rates. The experiments in \cref{sec:cfl} study the dependence of the CFL condition on the discretization parameters, while \cref{sec:td} highlight the time-domain convergence. In \cref{sec:eff} and \cref{sec:rob} we study the computational efficiency and robustness of the mass lumping approach. Finally we present in \cref{sec:ringres} a more practical example where we also apply perfectly matched layers to simulate an open domain.

\subsection{Eigenvalue problem}
\label{sec:disc_evp}
Since our main concern is the spatial discretization we perform experiments for the scalar discrete eigenvalue problem (EVP)  \eqref{eq:disc_evp} {with $\mu=\varepsilon\equiv 1$}.

This EVP is a discretization of the problem to find $\lambda\in\R$ and nontrivial $ H$ such that
\begin{align}
  \label{eq:scal_evp}
  -\Delta  H &= \lambda  H, &\text{on }&\Omega,\\
  \label{eq:scal_evp_bc}
   H &= 0,&\text{on }&\partial\Omega. 
\end{align}
For our experiments we choose the square $\Omega = [0,\pi]^2$ where eigenpairs are given by 
\begin{align}
   H_{n,k}&=\sin(nx)\sin(ky),&\lambda_{n,k}&=n^2+k^2,&n,k\in\N^+.
\end{align}
Equivalenty we could derive the corresponding system for the acoustic case with $\rho_0=c\equiv 1$ resulting in
\begin{align*}
  -\Delta  Q &= \lambda  Q, &\text{on }&\Omega,\\
   Q &= 0,&\text{on }&\partial\Omega. 
\end{align*}
with the same eigenpair solutions and discrete matrix eigenvalue problem. Thus the numerical experiments for the EVP corresponding to the acoustic initial boundary value problem exhibit the same behavior as the ones corresponding to the electromagnetic one studied in the following subsections.

\begin{figure}
  \begin{subfigure}{0.5\textwidth}
  \includegraphics[width=\textwidth]{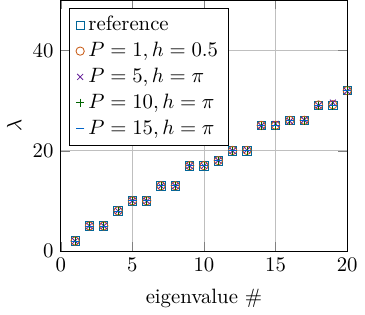}
  \caption{coarse mesh, high order}
  \end{subfigure}
  \begin{subfigure}{0.5\textwidth}
  \includegraphics[width=\textwidth]{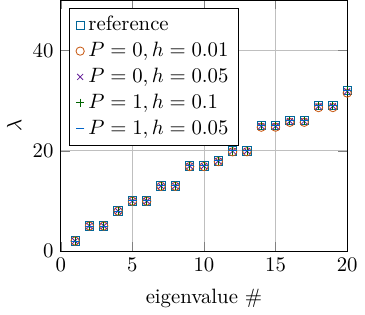}
  \caption{fine mesh, low order}
  \end{subfigure}
  \caption{Spectra of the discrete Dirichlet Laplacian on the square $[0,\pi]^2$ with polynomial degree $P$ and mesh size $h$.}
  \label{fig:disc_spec}
\end{figure}
\cref{fig:disc_spec} shows the discrete spectra of the different discretizations of the eigenvalue problem \cref{eq:scal_evp,eq:scal_evp_bc}. We observe no spurious eigenvalues for large polynomials degrees or finer mesh sizes respectively, i.e., all discrete eigenvalues converge towards their continuous counterparts.
\begin{figure}
  \begin{subfigure}{0.5\textwidth}
  \includegraphics[width=\textwidth]{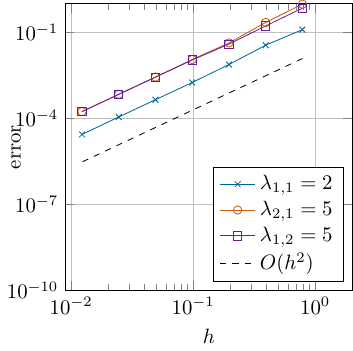}
  \caption{$ P=0$}
  \end{subfigure}
  \begin{subfigure}{0.5\textwidth}
  \includegraphics[width=\textwidth]{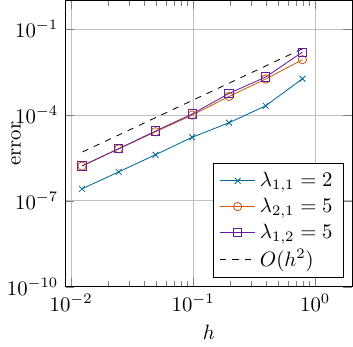}
  \caption{$ P=1$}
  \end{subfigure}\\
  \begin{subfigure}{0.5\textwidth}
  \includegraphics[width=\textwidth]{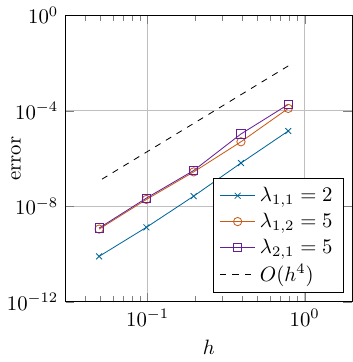}
  \caption{$ P=2$}
  \end{subfigure}
  \begin{subfigure}{0.5\textwidth}
  \includegraphics[width=\textwidth]{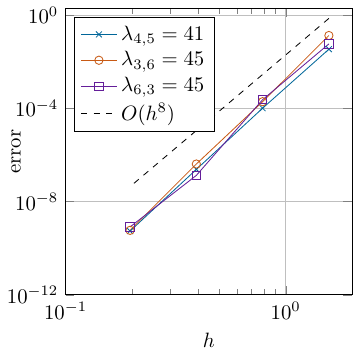}
  \caption{$ P=4$}
  \end{subfigure}
  \caption{Convergence of discrete eigenvalues for polynomial degrees $P$ and mesh sizes $h$.}
  \label{fig:ev_conv}
\end{figure}
\cref{fig:ev_conv} shows the convergence of selected discrete eigenvalues to their continuous counterparts. Note that while the experiments for degrees $P>0$ suggest a convergence rate of order $h^{2{P}}$ for the lowest order method we observe super convergence of the same order as for the method of order $1$ (i.e., a rate of $h^2$, which was already observed in \citep{kapidaniArbitraryorderCellMethod2021}). However the constant for the first order method improves significantly. Note that due to the high order convergence, in order to generate meaningful results (without hitting machine precision after very few refinements) we had to pick higher frequency resonances for the fourth order method.
\begin{figure}
  \includegraphics[width=\textwidth]{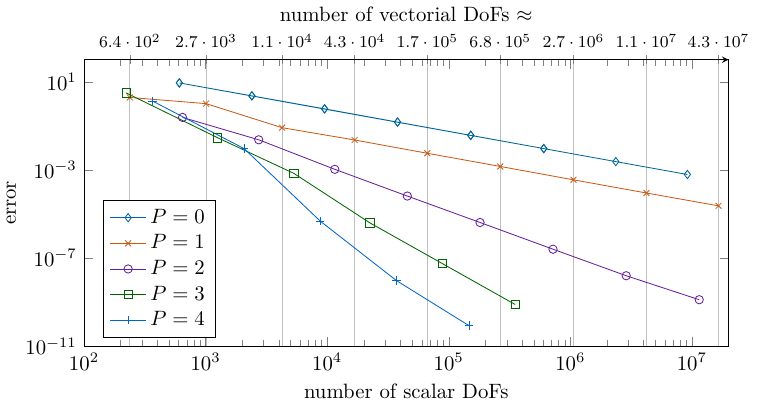}
  \caption{Convergence with respect to the primal number of DoFs, for the eigenvalue $\lambda_{8,3}=73$.}
  \label{fig:ev_conv_dof}
\end{figure}

\cref{fig:ev_conv_dof} shows the convergence with respect to the degrees of freedom for different polynomial degrees, suggesting that for smooth solutions a high order method is efficient. Also it confirms that, although we have quadratic convergence in $h$ for the polynomial degrees $ P=0,1$, the method with $P=1$ yields better results for the same number of degrees of freedoms.

\subsection{CFL condition}
\label{sec:cfl}
\begin{figure}
  \begin{subfigure}{0.5\textwidth}
    \includegraphics[width=\textwidth]{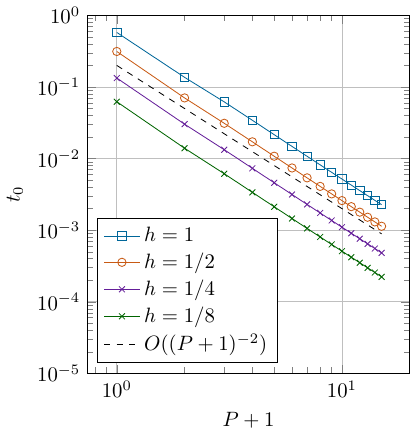}
  \caption{fixed mesh-size $h$}
  \end{subfigure}
  \begin{subfigure}{0.5\textwidth}
  \includegraphics[width=\textwidth]{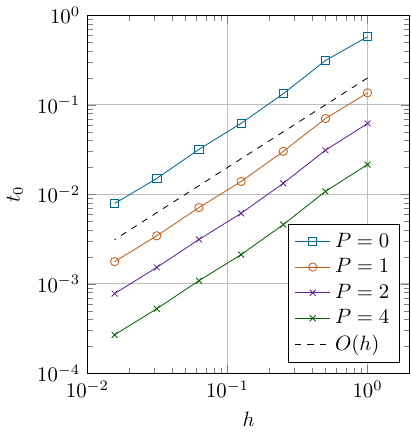}
  \caption{fixed polynomial degree $P$}
  \end{subfigure}
  \caption{The dependence of the maximal stable timestep $t_0$ on the discretization parameters.}
  \label{fig:cfl}
\end{figure}

Before we turn to time domain experiments we study numerically how the stability condition \cref{eq:cfl} on the timestep $\Delta t$ depends on the discretization parameters $h$ and $P$.
Figure \cref{fig:cfl} shows the dependence of the maximal stable timestep $t_0$ from \cref{eq:cfl} on the discretization parameters $P$ and $h$. We clearly observe that $t_0=O(h(P+1)^{-2})$ which is similar to the CFL condition for DG methods (cf., \citep[Chapter 4.7]{hesthavenNodalDiscontinuousGalerkin2008}).

\subsection{Time domain convergence}
\label{sec:td}
\begin{figure}
  \includegraphics[width=\textwidth]{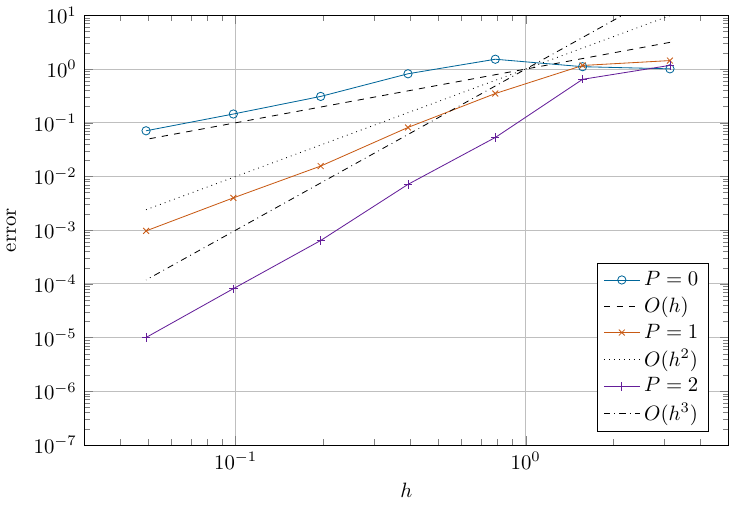}
  \caption{Convergence of time domain solutions with respect to the $\mathrm{L}^2$-norm.}
  \label{fig:td_conv}
\end{figure}

\cref{fig:td_conv} shows the convergence of the time domain solution for a problem with initial conditions 
\begin{align*}
   H_0(x,y)&=\sin(2x)\sin(6y),&\mathbf E_0(x,y)=0,
\end{align*}
at end time $T=1$ and a leap frog timestepping with time step sizes chosen such that the CFL condition is met (i.e., the discretization is stable) and the error is dominated by the spatial discretization. Other than for the eigenvalues we do not observe the super convergence for the lowest order case, but merely the expected rates of $h^{ P}$. 

\subsection{Efficiency}
\label{sec:eff}
\begin{figure}
  \includegraphics[width=\textwidth]{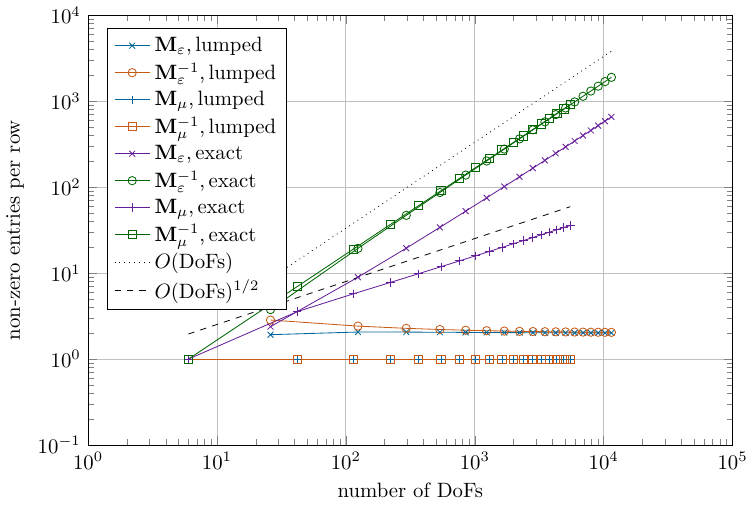}
  \caption{Number of non-zero entries per row of the (inverse) mass matrices obtained by using exact and lumped integration for a mesh with $h=0.5$ (six elements) and polynomial degrees $P$ from $0$ to $17$.}
  \label{fig:sparsity}
\end{figure}

To study the efficiency of using the mass lumped matrices we compare the number of non-zero entries of the mass matrices and the inverse matrices in the lumped and exact case respectively. To this end we choose a coarse mesh with six elements on the unit square and compare the number of non-zero entries in \cref{fig:sparsity} for varying polynomial degree.
We compare the matrices obtained by using the lumped inner products (cf., \cref{sec:inner_products}) to the ones obtained by exactly evaluating the integrals of the $\mathrm{L}^2$ inner products.
Since the lumped mass matrices are (block) diagonal with block sizes independent of the polynomial degree, we observe that the number of non-zero entries of the lumped mass matrices grows linearly with the number of degrees of freedom, even as the polynomial degree increases. The exact mass matrices are also block-diagonal, however, their block sizes grow with the polynomial degree (similar to the ones in \citep{kapidaniArbitraryorderCellMethod2021}). For the exact mass matrix of the scalar variable we observe a growth of the non-zero entries per row of power $1/2$. For the remaining (inverse) mass matrices we observe a linear growth of non-zero entries per row with respect to the number of unknowns.

\subsection{Robustness of computational efficiency}
\label{sec:rob}
\begin{figure}
  \begin{subfigure}{0.5\textwidth}
    \includegraphics[width=0.9\textwidth,clip,trim={0 150 0 150}]{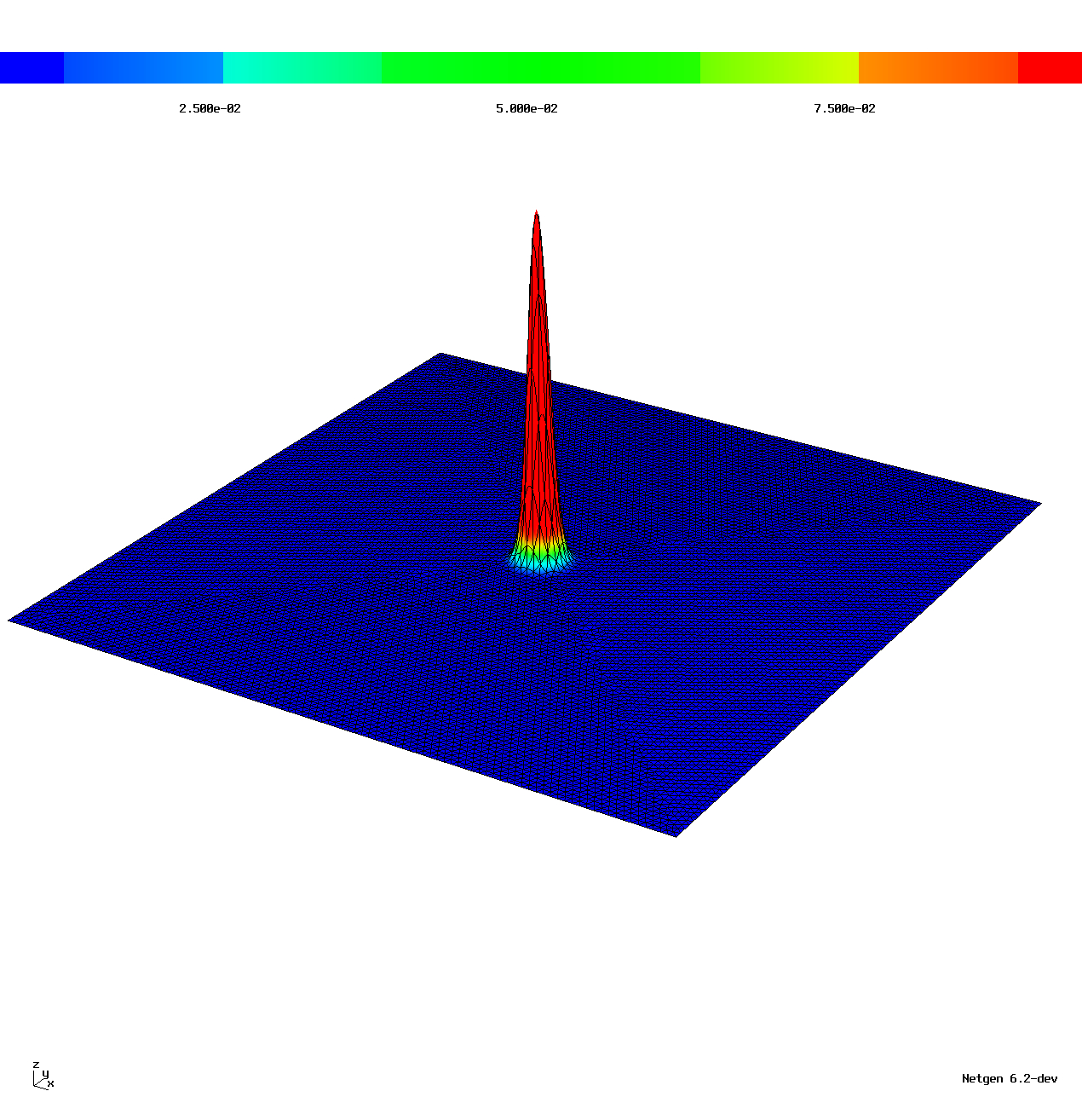}
    \caption{$T=0.0$}
  \end{subfigure}
%
  \begin{subfigure}{0.5\textwidth}
    \includegraphics[width=0.9\textwidth,clip,trim={0 150 0 150}]{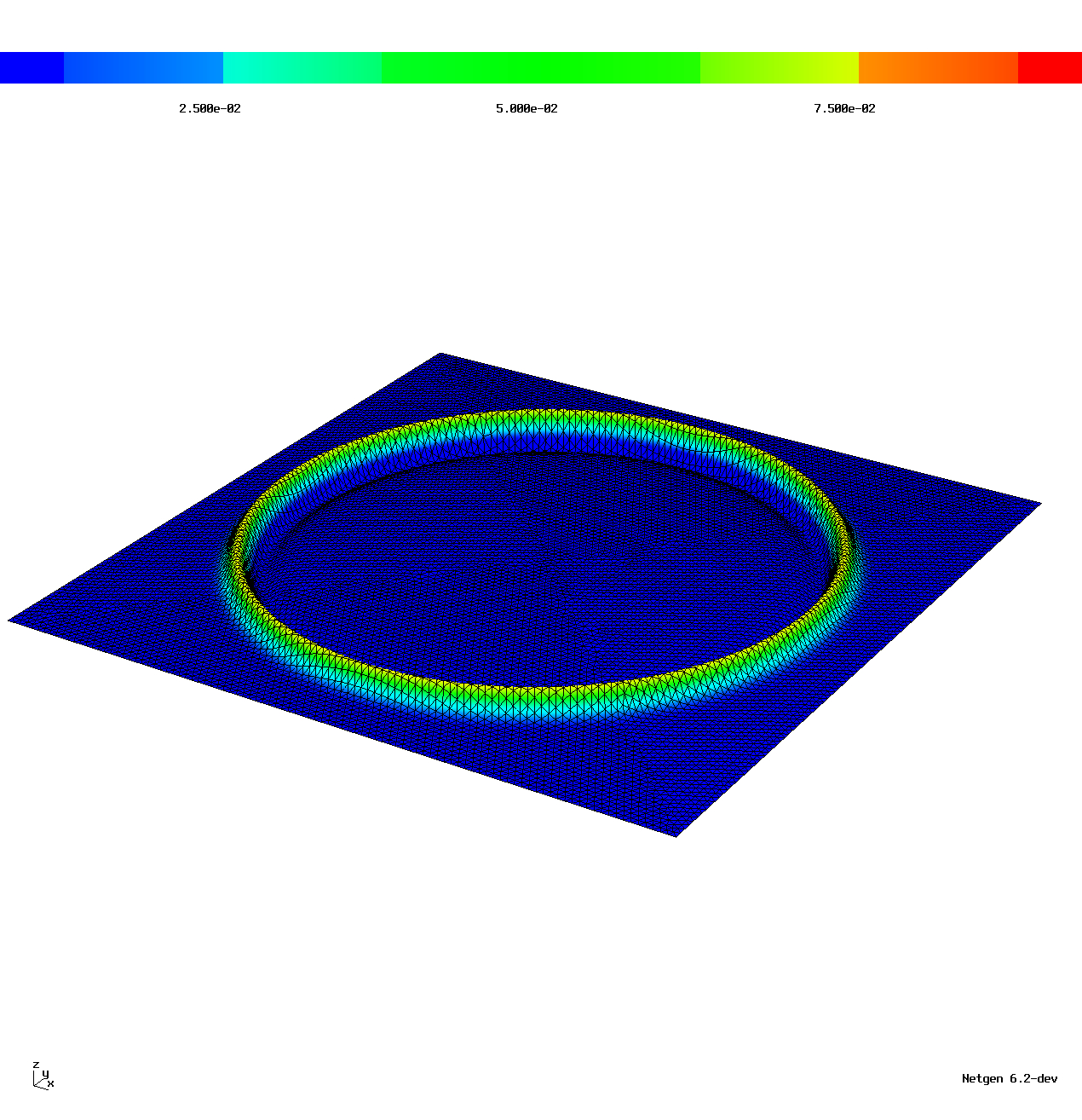}
    \caption{$T=0.4$}
  \end{subfigure}

  \begin{subfigure}{0.5\textwidth}
    \includegraphics[width=0.9\textwidth,clip,trim={0 150 0 150}]{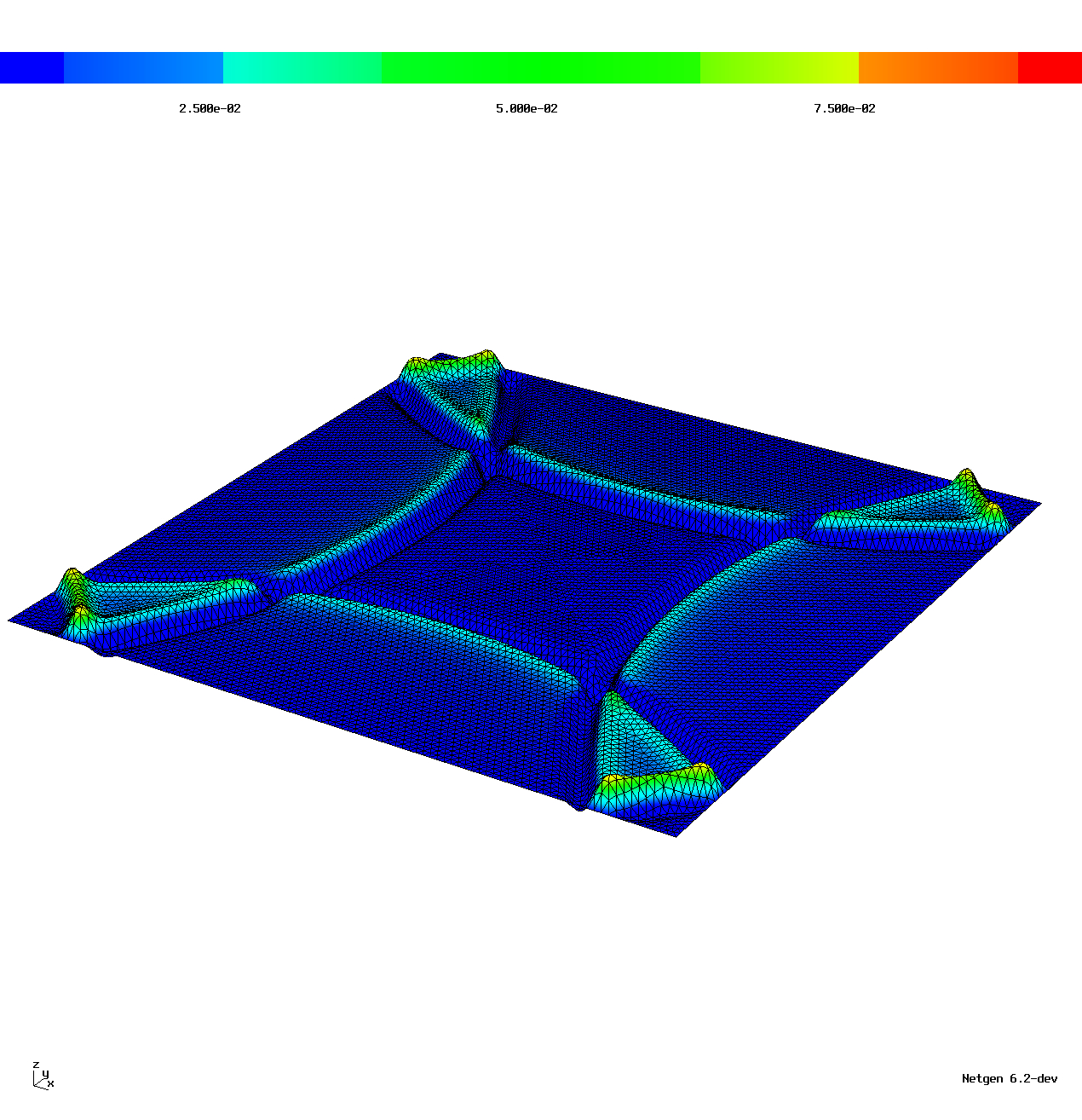}
    \caption{$T=0.8$}
  \end{subfigure}
%
  \begin{subfigure}{0.5\textwidth}
    \includegraphics[width=0.9\textwidth,clip,trim={0 150 0 150}]{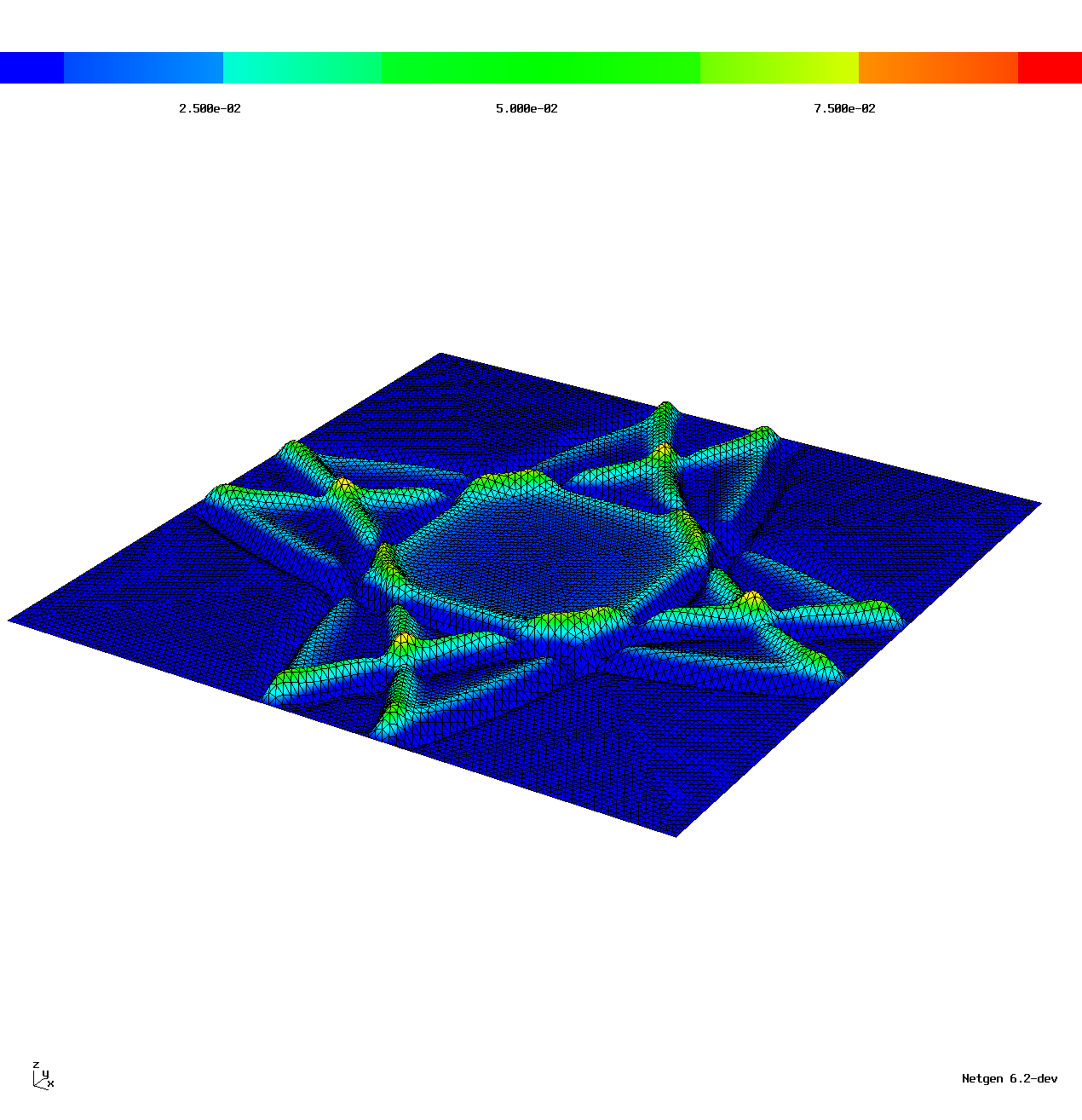}
    \caption{$T=1.2$}
  \end{subfigure}

  \begin{subfigure}{0.5\textwidth}
    \includegraphics[width=0.9\textwidth,clip,trim={0 150 0 150}]{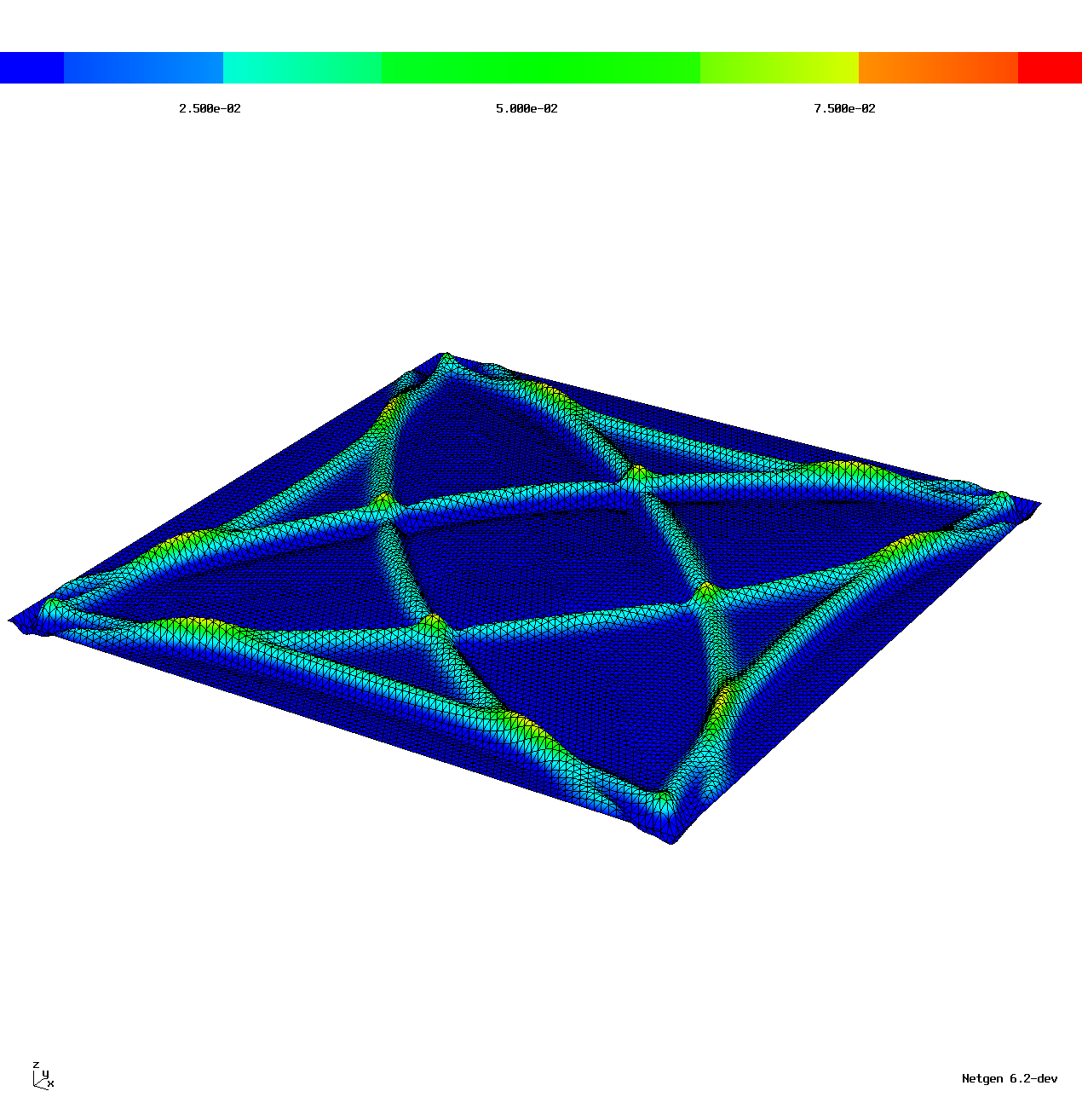}
    \caption{$T=1.6$}
  \end{subfigure}
%
  \begin{subfigure}{0.5\textwidth}
    \includegraphics[width=0.9\textwidth,clip,trim={0 150 0 150}]{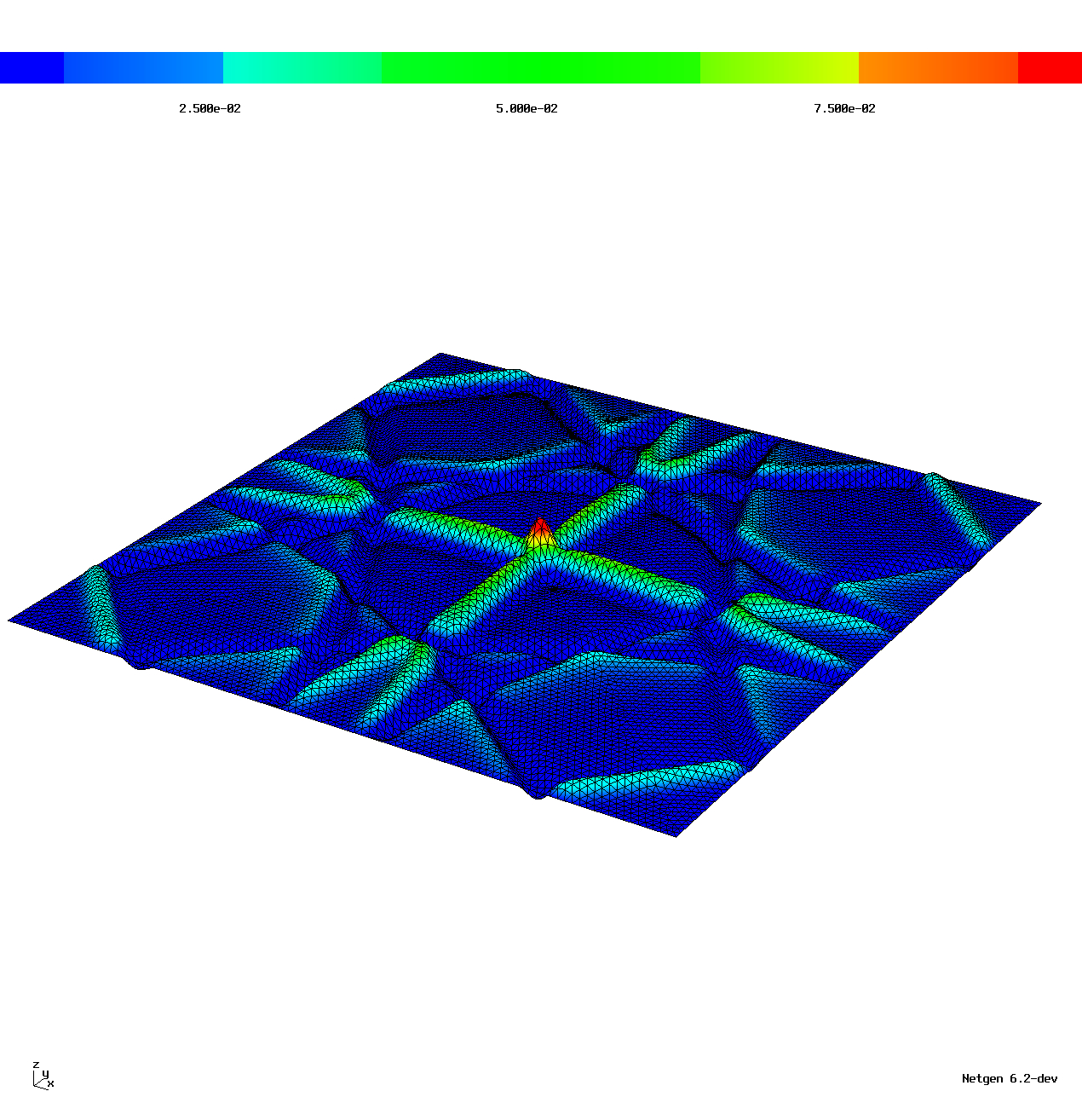}
    \caption{$T=2.0$}
  \end{subfigure}
  \caption{Snapshots of the time domain solution for an inital Gauss peak for $t\in[0,2]$.}
  \label{fig:gauss_peak}
\end{figure}

Due to the matrix-free nature of our method we expect to be able to do large scale computations without significant memory requirements\footnote{All computations were carried out on an off-the shelf desktop computer with 4 CPUs with 3.3GHz and 16GiB of Memory.}. 
Thus we expect the number of unknowns which can be computed per second to be independent of the number of finite elements (i.e., the mesh-size for a given problem). The efficiency of the inverse mass matrices is also independent of the polynomial degree while applying the discrete differential operators has a (mild) dependence on the polynomial degree.
To test this we run a sequence of time dependent problems on the unit square with initial fields
\begin{align*}
  H_0(x,y) &= \exp\left(-50^2\left((x-0.5)^2+(y-0.5)^2\right)\right),&
  \mathbf E_0(x,y) &= 0,
\end{align*}
for different mesh-sizes and polynomial degrees. Snapshots of the resulting solution for $h=0.01$, $ P=4$ is shown in \cref{fig:gauss_peak}. The step-sizes $\Delta t$ are chosen experimentally such that the resulting discretization is stable.
Table \ref{tab:dofspersec} gives the according numbers of DoFs per second, where for each entry the minimal time of computing four times 50 steps is given. While we observe a small decay of the number of DoFs/s for larger problems (probably due to the more costly application of the differential operators), the results confirm that for the given problem sizes the number of DoFs/s is roughly $10^8$, mostly independent of the mesh-size $h$ and polynomial $P$ degree.

On a more qualitative note it can be observed in \cref{fig:gauss_peak} that the width of the transported peak is preserved during the computation, i.e., we do not observe any numerical dispersion.

\begin{table}
  \pgfplotstabletypeset[
    columns = {maxh,order,tau,scalardofs,vectorialdofs,totaldofs,dofspers},
    columns/maxh/.style={column name =$h$},
    columns/order/.style={column name =$P$,int detect,},
    columns/tau/.style={column name =$\Delta t$,sci},
    columns/scalardofs/.style={column name =scal. DoFs,sci},
    columns/vectorialdofs/.style={column name =vect. DoFs,sci},
    columns/totaldofs/.style={column name =total DoFs,sci},
    columns/dofspers/.style={column name =DoFs/s,sci,postproc cell content/.style={
      @cell content/.add={$\bf}{$}}},
    column type/.add={|}{},
    every head row/.style={after row=\hline},
  ]{data/peak_dofspersecond_50_4.out}
  \caption{Computed DoFs/s for the initial value problem for different discretizations.}
  \label{tab:dofspersec}
\end{table}

\subsection{Ring resonator}
\label{sec:ringres}
\begin{figure}
  \begin{subfigure}{0.5\textwidth}
    \includegraphics[width=\textwidth]{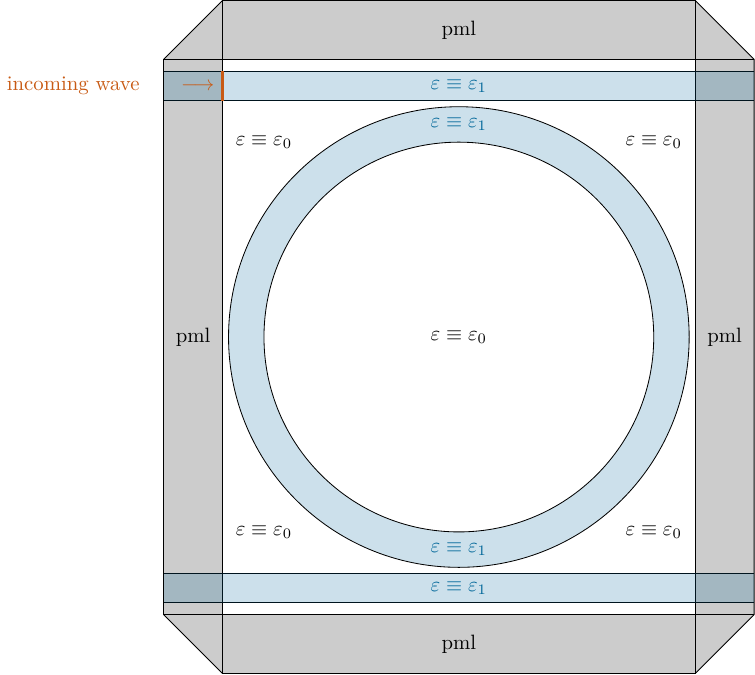}
  \caption{geometry}
    \label{fig:ring_resonator_geo}
  \end{subfigure}
  \begin{subfigure}{0.5\textwidth}
    \includegraphics[width=\textwidth]{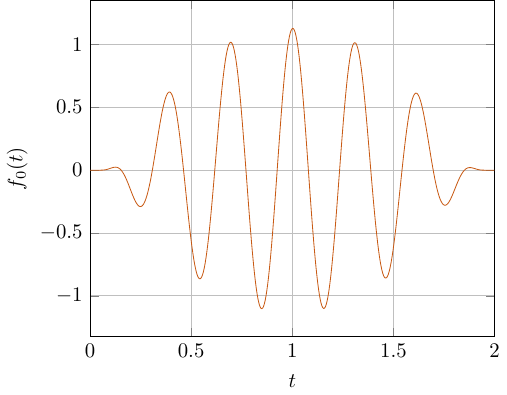}
  \caption{source}
    \label{fig:ring_resonator_source}
  \end{subfigure}
  \caption{Configuration of the ring resonator experiment.}
  \label{fig:ring_resonator}
\end{figure}

To test our method in a more challenging setting we choose a model of a ring resonator, which is given by two parallel electrical wires with a looped wire in between (cf. \cref{fig:ring_resonator_geo}). A wave is introduced at the left end of the top wire (cf. \cref{fig:ring_resonator_source}) which resonates in the loop and induces an output at the left end of the bottom wire.

The unbounded physical domain is modelled by the use of a perfectly matched layer (pml) (cf. \citep{berengerPerfectlyMatchedLayer1994,teixeiraComplexSpaceApproach2000}).

\begin{figure}
  \begin{subfigure}{0.5\textwidth}
    \includegraphics[width=\textwidth,clip,trim = {20 100 20 200}]{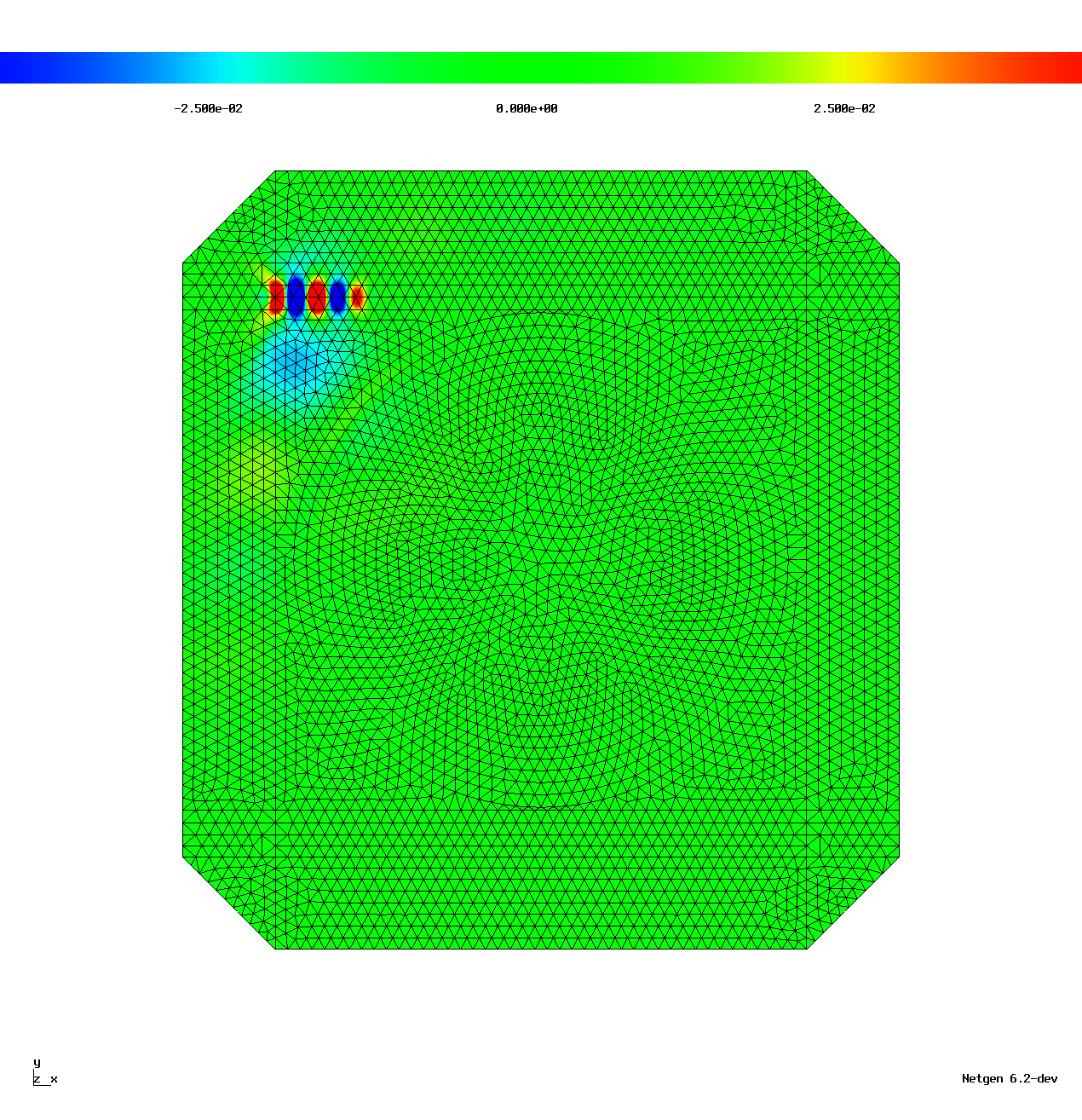}
  \caption{$t=1$}
  \end{subfigure}
  \begin{subfigure}{0.5\textwidth}
    \includegraphics[width=\textwidth,clip,trim = {20 100 20 200}]{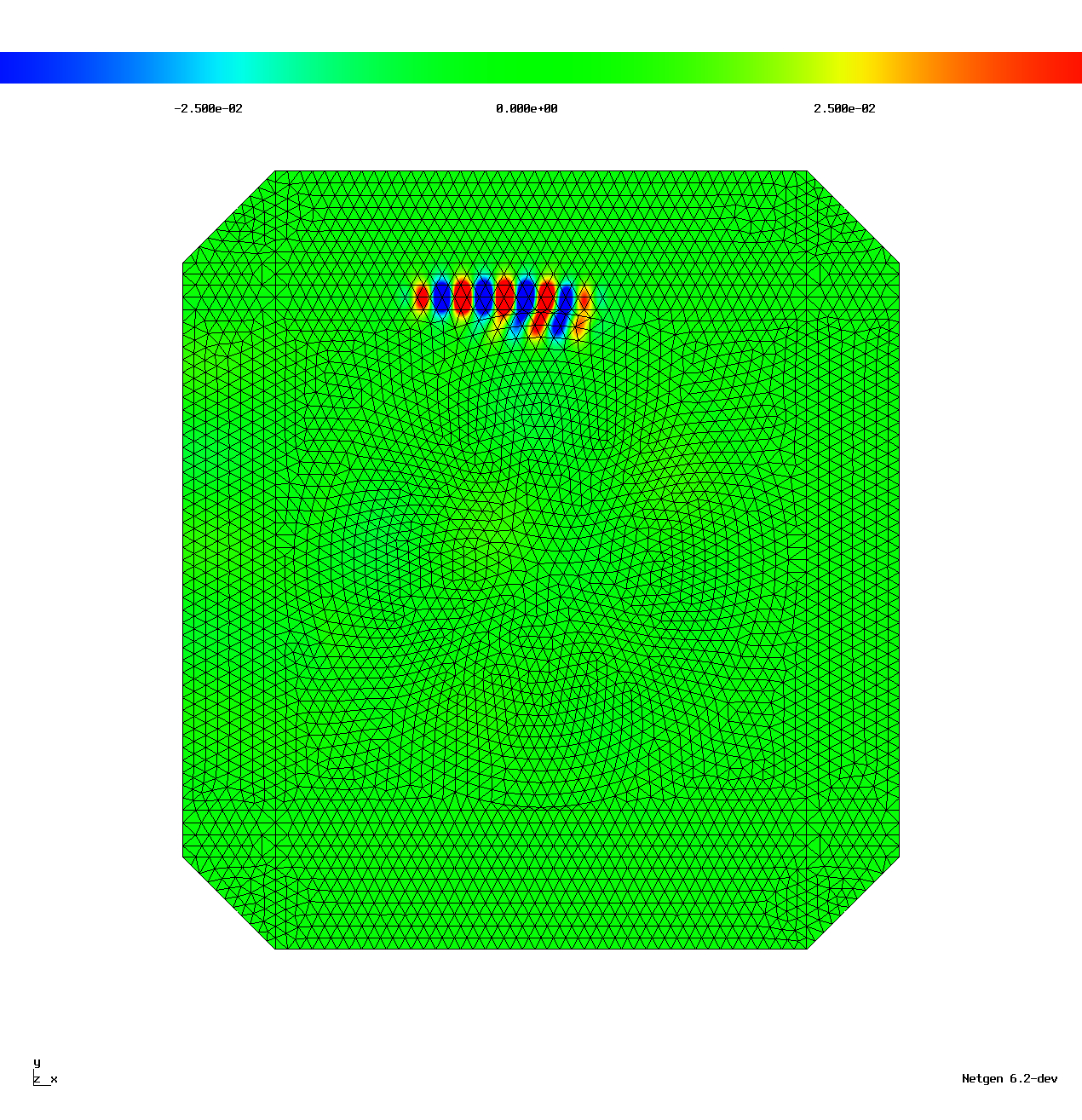}
  \caption{$t=3$}
  \end{subfigure}
  \begin{subfigure}{0.5\textwidth}
    \includegraphics[width=\textwidth,clip,trim = {20 100 20 200}]{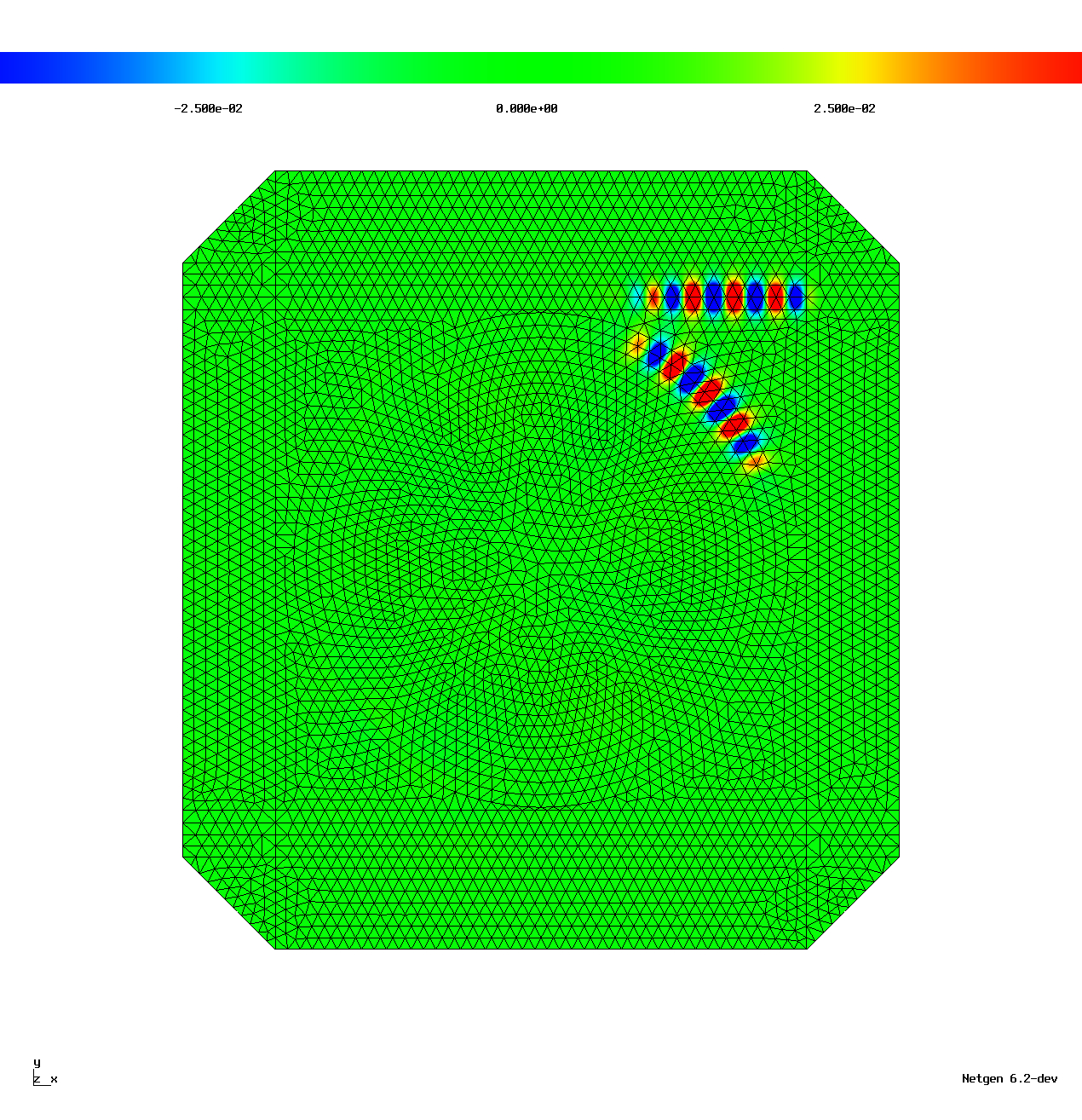}
  \caption{$t=5$}
  \end{subfigure}
  \begin{subfigure}{0.5\textwidth}
    \includegraphics[width=\textwidth,clip,trim = {20 100 20 200}]{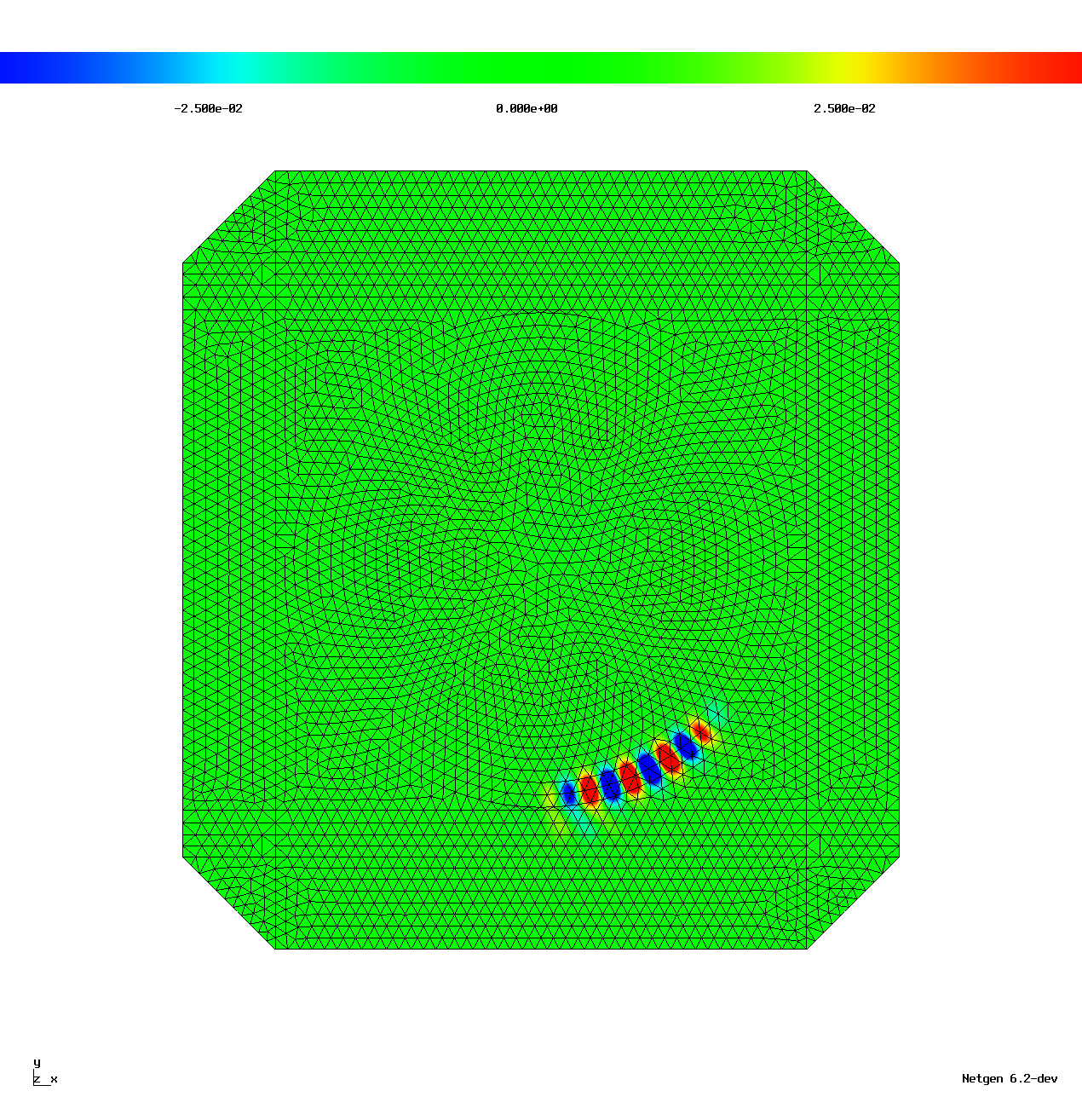}
  \caption{$t=9$}
  \end{subfigure}
  \begin{subfigure}{0.5\textwidth}
    \includegraphics[width=\textwidth,clip,trim = {20 100 20 200}]{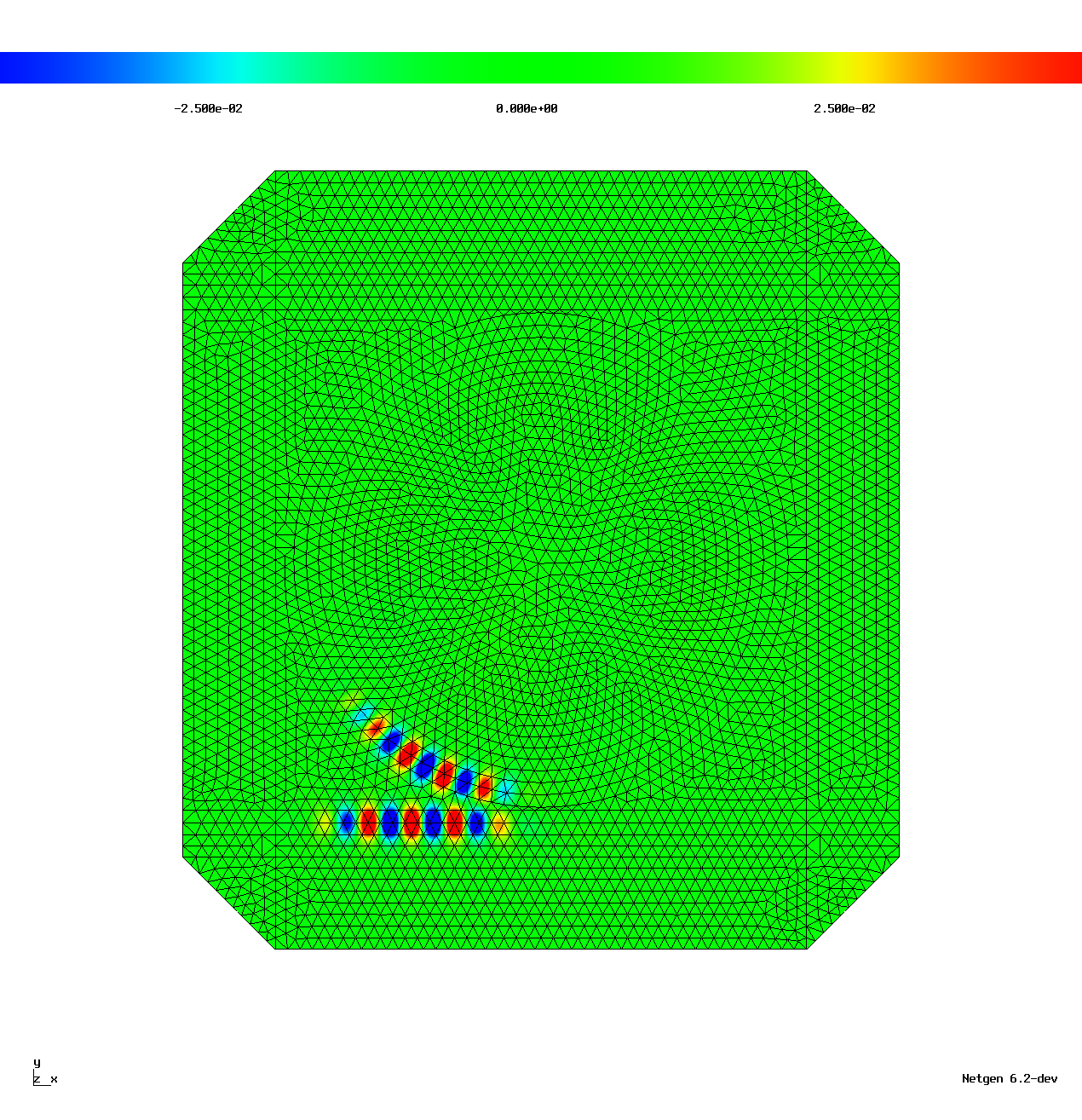}
  \caption{$t=11$}
  \end{subfigure}
  \begin{subfigure}{0.5\textwidth}
    \includegraphics[width=\textwidth,clip,trim = {20 100 20 200}]{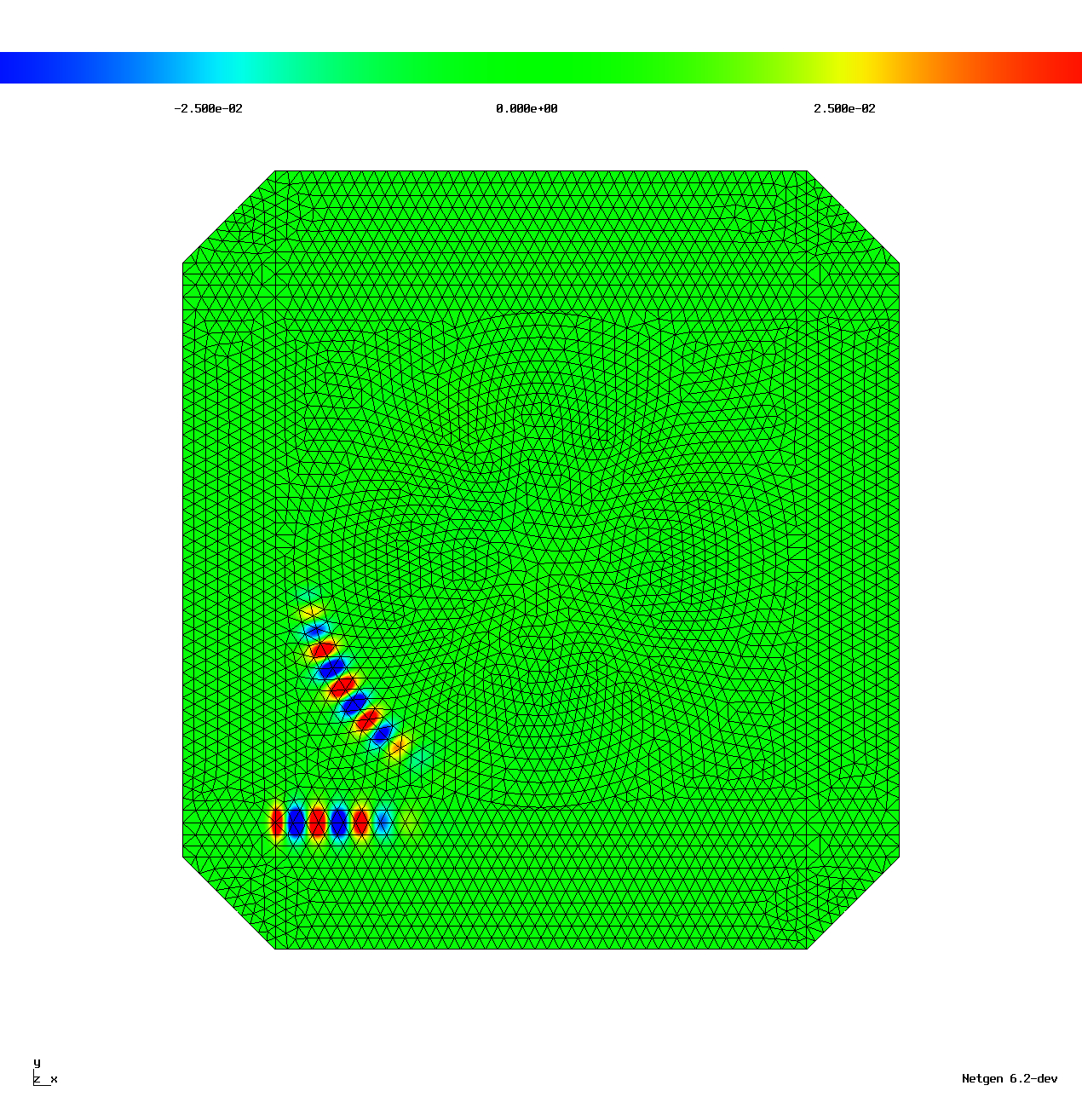}
  \caption{$t=12$}
  \end{subfigure}
  \caption{Snapshots of the ring resonator experiments (cf., \cref{fig:ring_resonator}).}
  \label{fig:ring_resonator_res}
\end{figure}

To obtain the results shown in \cref{fig:ring_resonator_res}  we use a mesh with 
$\approx 7.8\cdot 10^3$ finite elements which results in spaces of dimension 
$\approx 7.1\cdot 10^5$ for the scalar and  
$\approx 1.5\cdot 10^6$ for the vectorial space when using a basis of order $ P=5$. Note that due to the auxiliary unknowns of the PML formulation the number of total unknowns is even higher. We choose a time-step size $\Delta t=1.15\cdot 10^{-4}$ which is experimentally confirmed to be stable (i.e., the CFL condition is fulfilled). Thus to compute up to $t=12$, 
a number of $\approx 1.0\cdot 10^5$ timesteps have to be computed. All computations where carried out on a desktop computer in less than a day of computation time, with one timestep taking $\approx 0.46s$. Due to the efficiency of the mass lumped basis the amount of used memory and time to set up the problem is negligable since no large (inverse) matrix has to be set up/stored.

\section{Conclusion and Outlook}\label{sec:conclusion}
  In the work at hand we have presented a new approach to construct basis functions for an arbitrary order accurate (in the {spacial} discretisation) cell method. In particular, this new basis resolves some stability issues from previous approaches and has the additional property that mass lumping is easily available for the respective mass matrices. Therefore our new approach results in well-conditioned (inverse) mass matrices {with} a very favourable sparsity pattern which is uniform in the polynomial degree of the basis.
  We {introduced} the resulting mass lumped formulation both, for the time domain Maxwell system and for the acoustic wave equation.
We have applied our method in numerous numerical experiments to underline these claims and shown that the method scales very well for larger problems even on off-the-shelf desktop computers.
  Future work includes the extension of our method to higher dimensions for acoustic, as well as electromagnetic problems. We aim at optimising the simulation of unbounded domains by using more {sophisticated} PML constructions such as the ones based on infinite elements of \citep{bermudezOptimalPerfectlyMatched2007,nannenComplexscaledInfiniteElements2022} and exploring alternative approximation techniques on quadrilateral elements such as spline approximations of curved boundaries, e.g., \citep{ratnaniArbitraryHighOrderSpline2012} or \citep{kapidaniHighOrderGeometric2023}. 
Finally, the theoretical numerical analysis of the method is being carried out and will be presented elsewhere.

\section*{Acknowledgements}
Author Bernard Kapidani has been partially supported by the Swiss National Science Foundation via the project HOGAEMS n.200021\_188589.

%
%
%
\bibliographystyle{unsrt}
\bibliography{biblio}

\end{document}